\documentclass[a4paper,12pt]{article}

\scrollmode
\usepackage{amsmath,amssymb,amsfonts,graphicx,psfrag,subfig}
\usepackage[latin1]{inputenc}

\newtheorem{thm}{Theorem}[section]
\newtheorem{prop}[thm]{Proposition}
\newtheorem{lem}[thm]{Lemma}
\newtheorem{df}[thm]{Definition}

\newtheorem{rem}[thm]{Remark}

\newtheorem{cor}[thm]{Corollary}

\def\be#1 {\begin{equation} \label{#1}}
\newcommand{\ee}{\end{equation}}
\def\dem {\noindent {\bf Proof: }}

\newcommand{\mb}{\medskip\noindent}
\newcommand{\gb}{\bigskip\noindent}
\newcommand{\R}{\mathbb R}

\newcommand{\Z}{\mathbb Z}

\newcommand{\B}{\mathcal B}
\newcommand{\F}{\mathcal F}

\newcommand{\ww}{\langle}

\def \NNN {\mathcal{N}}
\def \CCC {\mathcal{C}}
\def \eee {\mathrm{e}}

\def \qqq {{\bf q}}
\def \qq {\mathrm{q}}

\def \UU {{\bf U}}

\def \GG {{\bf G}}

\def \ww {{\bf w}}
\def \vv {{\bf v}}

\def \NN {\mathrm{N}}

\def \virg {\, , \,\,}
\def \dsp {\displaystyle}
\def \vsp {\vspace{6pt}}

\def\sqw{\hbox{\rlap{\leavevmode\raise.3ex\hbox{$\sqcap$}}$%
\sqcup$}}
\def\findem{\ifmmode\sqw\else{\ifhmode\unskip\fi\nobreak\hfil
\penalty50\hskip1em\null\nobreak\hfil\sqw
\parfillskip=0pt\finalhyphendemerits=0\endgraf}\fi}

\title{Existence of sweeping process in Banach spaces under directional prox-regularity.}
\date{\today }
\author{ Fr\'ed\'eric Bernicot\\ Universit\'e de Paris-Sud\\F-91405 Orsay
Cedex\\frederic.bernicot@math.u-psud.fr \and Juliette Venel\\ Universit\'e de Paris-Sud\\F-91405 Orsay
Cedex\\juliette.venel@math.u-psud.fr}

\begin {document}

\maketitle

\begin{abstract} This paper is devoted to weaken ``classical'' assumptions and to give new arguments to prove existence of sweeping process (associated to the proximal normal cone of sets). Mainly we define the concept of ``directional prox-regularity'' and give assumptions on a Banach space to ensure the existence of such sweeping process (which permits to generalize existing results requiring a Hilbertian structure). 
\end{abstract}

\mb {\bf Key-words:} Sweeping process - Prox-regularity - Proximal Normal cone. 

\mb {\bf MSC:} 34A60, 49J52, 58C20.

\tableofcontents

\gb

\section{Introduction}

The aim of this paper is to prove existence results for sweeping
process associated to a moving set $t\to C(t)$ on a time interval
$I:=[0,T]$ in considering the ``proximal normal cone'' (which will be
denoted $\NN$). Let $\B$ be a Banach space and $C:I \rightrightarrows
\B$ be a set-valued map with nonempty closed values, and let $F: I \times \B \rightrightarrows \B$ be a set-valued map taking nonempty convex compact values. An associated sweeping process $u:I\to \B$ is a solution of the following differential inclusion:
\begin{equation} 
\left\{
\begin{array}{l}
 \dsp \frac{du(t)}{dt} + \NN(C(t),u(t)) \supset F(t,u(t)) \vsp \\
 u(t)\in C(t) \vsp \\
 u(0)=u_0\ ,
\end{array}
\right. \label{sys1}
\end{equation}
with an initial data $u_0\in C(0)$. This differential inclusion can be thought as following: the point $u(t)$, submitted to the field $F(t,u(t))$, has to live in the set $C(t)$ and so follows its time-evolution.

\gb We begin first by detailing the story of the study for similar problems. The sweeping processes have been introduced by
 J.J. Moreau in 70's  (see \cite{Moreausweep}). He considered the following problem: a point $u(t)$ has to be inside a moving convex set $C(t)$ included in a Hilbert space. When this point is catched-up by the boundary of $C(t)$, it moves in the opposite of the outward normal direction of the boundary, as it is pushed by the physical boundary in order to stay inside the convex set $C(t)$. Then the position $u(t)$ of this point is described by the following differential inclusion
\begin{equation}
\label{Moreq}
 -\dot{u}(t) \in \partial I_{C(t)}(u(t)).
\end{equation}
Here we write $\partial I_{C}$ for the subdifferential of the indicator function of a convex set $C$. In this work, the sets $C(t)$ are assumed to be convex and so $\partial I_{C(t)}$ is a maximal monotone operator depending on time.
To solve this problem, J.J. Moreau brings a new important idea in proposing a \textit{catching-up} algorithm.
To prove the existence of solutions, he builds discretized solutions in dividing the time interval $I$ into sub-intervals where the convex set $C$ does not vary too much. Then by compactness arguments, he shows that we can construct a limit mapping (when the length of subintervals tends to $0$) which satisfies the desired differential inclusion. \\
Indeed with well-known convex analysis, as $C(t)$ is convex, we have
$\partial I_{C(t)}(x) = \NN(C(t),x)$. So it is the first result
concerning sweeping process (with no perturbation $F=0$).

\gb Since then, important improvements have been developped by weakening the assumptions in order to
obtain the most general result of existence for sweeping
process. There are several directions: one can want to add a perturbation
$F$ as written in (\ref{sys1}), one may require a weaker assumption than the convexity of the sets, one would like to obtain results in Banach spaces (and not only in Hilbert spaces), ...

\gb In \cite{V}, M.~Valadier dealt with sweeping process by sets $C(t)=\R^n\setminus \textrm{int}(K(t))$ where $K(t)$ are closed and convex sets.
Then in \cite{CXDV}, C.~Castaing, T.X. D\'uc H\={a} and M. Valadier have studied the perturbed problem in finite dimension ($\B=\R^n$) with convex sets $C(t)$ (or complements of convex sets). In this framework, they proved existence of solutions for (\ref{sys1}) with a convex compact valued perturbation $F$ and a Lipschitzean multifunction $C$. 
Later in \cite{Castaingper}, C.~Castaing and M.D.P.~Monteiro Marques have considered similar problems in assuming upper semicontinuity for $F$ and a ``linear compact growth'':
\begin{equation}
 F(t,x)\subset \beta(t) (1+|x|)\overline{B(0,1)} \virg \forall (t,x) \in I \times \R^n.
\label{hypF}
\end{equation}
Moreover the set-valued map $C$ was supposed to be Hausdorff continuous and satisfying an ``interior ball condition'': 
\begin{equation} 
 \exists r>0 \virg  B(0,r) \subset C(t) \virg \forall t \in I.
\label{hypC}
\end{equation}

\gb  Then the main concept, which appeared to get around the convexity of sets $C(t)$, is the notion of  ``uniform prox-regularity''. This property is very well-adapted to the resolution of (\ref{sys1}): a set $C$ is said to be {\it $\eta$-prox-regular} if the projection on $C$ is single valued and continuous at any point whose the distance to $C$ is smaller than $\eta$. \\
Numerous works have been devoted to applications of prox-regularity in the study of sweeping process. The case without perturbation ($F={0}$) was firstly treated by G.~Colombo, V.V.~Goncharov in \cite{Colombo}, by H.~Benabdellah in \cite{Benab} and later by L.~Thibault in \cite{Thibsweep} and by G. Colombo, M.D.P.~Monteiro Marques in \cite{Monteiro}.
In \cite{Thibsweep}, the considered problem is 
\begin{equation}
 \left \{
\begin{array}{l} 
-du \in \NN(C(t), u(t)) \vsp \\
u(T_0)=u_0\ ,
 \end{array}
\right.
\label{eq:mesdiffssm}
\end{equation} 
where $du$ is the differential measure of $u$. The existence and uniqueness of solutions of (\ref{eq:mesdiffssm}) are proved with similar assumptions as previously. \\
In infinite dimension with assuming that $\B$ is a Hilbert space $\B=H$, the perturbed problem is studied by M.~Bounkhel, J.F.~Edmond and L.~Thibault in \cite{Thibnonconv, Thibsweep, Thibrelax, Thibbv} (see Theorem \ref{thm1}).
For example in \cite{Thibbv}, the authors show the well-posedness of
\begin{equation}
 \left \{
\begin{array}{l} 
-du \in \NN(C(t), u(t)) + F(t,u(t)) dt \vsp \\
 u(0)=x_0\ ,
 \end{array}
\right.
\label{eq:mesdiffasm}
\end{equation} 
with  a set-valued map $C$ taking $\eta$-prox regular values (for some $\eta>0$) such that   
\begin{equation}
|d_{C(t)}(y) - d_{C(s)}(y)| \leq \mu(]s,t]) \virg \forall y\in H,\ \forall \  s, t \in I  \virg s\leq t
\label{varadon}
\end{equation}
where $\mu$ is a nonnegative measure satisfying
\begin{equation}
\sup_{s \in I} \mu(\{s\}) <\frac{\eta}{2}.
 \label{charge_singleton}
\end{equation}
The proof use the algorithm developped by J.J.~Moreau with additional arguments to deal with the prox-regularity assumption.

\gb Indeed the main difficulty of this problem is the weak smoothness of the proximal normal cone. For a fixed closed subset $C$, the set-valued map $x\rightarrow \NN(C,x)$ is not upper semicontinuous, which is needed for the proof. The prox-regularity implies this required smoothness. We finish by presenting the work of H.~ Benabdellah (see \cite{bena2}). He deals with sweeping process in an abstract Banach framework, in considering the Clarke normal cone, which satisfies this upper semicontinuity. 

\gb After this description of existing results, we come to our contribution in this article. We are looking for results concerning sweeping process with proximal normal cone. We first  precise some results (essentially already obtained in the previously cited papers) about these ones in an Hilbert framework.
Then in Section \ref{sec:ex} we explain with an example due to a model of crowd motion
(detailed in the thesis of one of us \cite{TheseJu}) that the
``uniform prox-regularity'' assumption could fail for some interesting
cases. We define also in Subsection \ref{subsec:dir} a weaker notion, which corresponds to a ``directional
prox-regularity'' property. Moreover we present new arguments for the
proof of existence of sweeping process. It is still based on the
ideas of the catching-up algorithm of J.J. Moreau. This algorithm gives us a sequence of functions
(corresponding to discretized solutions), whose we can extract a
weak-convergent subsequence. The technical problem is to check that
this limit function is a solution of the differential inclusion. The well-known arguments use the Hilbertian structure of the space, and the fact that the support function of the subdifferential of the distance function is upper semicontinuous (which is implied by the prox-regularity of the set).
Here we propose a new approach to describe this ``weak continuity''. This allows us to present results in an abstract Banach framework (under some assumptions on the Banach space, see Subsection \ref{subsec:geo}) and to deal only with a ``directional prox-regularity''. We describe these new arguments for a single-valued perturbation $F$, which will be denoted by $f$. Here are our two main results (proved in Section \ref{sec:res}):

\begin{thm} \label{thm:general} 
Let $\B$ be a separable, reflexive, uniformly smooth Banach space, which is ``$I$-smoothly weakly compact'' for an exponent $p\in[2,\infty)$ (see Definition \ref{df:imp}).
Let $f: \B \rightarrow \B$ be a continuous function admitting at most a linear growth and $r>0$ be a fixed real. Let $C$ be a nonempty ball-compact $(r,f)$ prox-regular subset of $\B$. 
Then for all $u_0\in C$, the system 
$$ \left\{ \begin{array}{ll}
 \dot{u}(t)+\NN(C,u(t)) \ni f(u(t)) \\
 u(0)=u_0
\end{array} \right. $$
has an absolutely continuous solution $u$ and for all $t\in I$, $u(t)\in C$.
\end{thm}

\mb In the case of a Hilbert space $\B=H$, we do not need to require the ball-compactness of the set $C$ and prove~:

\begin{thm} \label{thm:genera} 
Let $\B=H$ be a separable Hilbert space. Let $f: \B \rightarrow \B$ be a Lipschitz function admitting at most a linear growth and $r>0$ be a fixed real. Let $C$ be a nonempty closed $(r,f)$ prox-regular subset of $H$. 
Then for all $u_0\in C$, the system 
$$ \left\{ \begin{array}{ll}
 \dot{u}(t)+\NN(C,u(t)) \ni f(u(t)) \\
 u(0)=u_0
\end{array} \right. $$
has one and only one absolutely continuous solution $u$ and for all $t\in I$, $u(t)\in C$.
\end{thm}

\mb Moreover we give in the last subsections several extensions in a Banach and Hilbert framework with a non constant set-valued map $C$. 

\section{Preliminaries}

For an easy reference, we recall the main definitions and notation, used throughout the paper. Let $\B$ be Banach space, equipped with its norm $\|\ \|$. We write $B(x,r)$ for the open ball of center $x$ and of radius $r$ and $\overline{B}(x,r)$ for its closure. Let $S_1 $ and $S_2 $ be two nonempty subsets of $\B$, we denote their Hausdorff distance $H(S_1,S_2) $ defined by $$H(S_1,S_2):= \max\left( \sup_{x \in S_1} d(x,S_2),\sup_{x \in S_2} d(x,S_1) \right)   .$$
\begin{df} Let $C$ be a closed subset of $\B$. The set-valued projection operator $P_C$ is defined by
$$ \forall x\in \B, \qquad P_C(x):=\left\{y\in C,\ \|x-y\|=d(x,C)\right\}.$$
\end{df}

\begin{df} Let $C$ be a closed subset of $\B$ and $x\in C$, we denote by $\NN(C,x)$ the proximal normal cone of $C$ at $x$, defined by:
 $$ \NN(C,x):=\left\{v\in\B,\ \exists s>0, \ x\in P_C(x+sv)\right\}.$$
\end{df}

\mb We now come to the main notion of ``prox-regularity''. It was
initially introduced by H. Federer (in \cite{Federer}) in spaces of finite dimension under the name of ``positively reached sets''. Then it was extended in Hilbert spaces by A. Canino in \cite{C} and A.S. Shapiro in \cite{Shapiro}. After, this notion was studied by F.H.~Clarke, R.J.~Stern and P.R.~Wolenski in \cite{Clarke} and by R.A. Poliquin, R.T.~Rockafellar and L.~Thibault in \cite{PRT}. Few years later, F. Bernard and L. Thibault have defined this notion in Banach spaces (see \cite{BT}).

\begin{df} Let $C$ be a closed subset of $\B$ and $r>0$. The set $C$ is said $\eta$-prox-regular if for all $x\in C$ and all $v\in \NN(C,x)\setminus \{0\}$
 $$B\left(x+\eta\frac{v}{\|v\|},\eta\right) \cap C = \emptyset.$$
\end{df}

\mb We refer the reader to \cite{Clarke,Clarke2} for other equivalent definitions related to the limiting normal cone. Moreover we can define this notion using the smoothness of the function distance $d(.,C)$, see \cite{PRT}.
This definition is very geometric, it describes the fact that we can
continuously roll an external ball of radius $\eta$ on the whole boundary of the set $C$. The main property is the following one: for an $\eta$-prox-regular set $C$, and for every $x$ satisfying $d(x,C)<\eta$, the projection of $x$ onto $C$ is well-defined and continuous.

\section{Some details about sweeping process in Hilbert spaces.} \label{sec1}

We consider the following ``sweeping process'' on a time interval $I=[0,T]$ with a single-valued perturbation $f$:

\be{eq:syst1}
\left\{ \begin{array}{ll}
 \dot{u}(t) + \NN(C(t),u(t)) \ni f(t,u),\qquad a.e.\ t\in I \vsp \\
 u(0)=u_0
\end{array}
 \right. 
\ee

\mb We recall the results of J.F. Edmond and L. Thibault (see Theorem 5.1 of \cite{Thibbv}):

\begin{thm} \label{thm1} Let $H$ be a Hilbert space, $\eta>0$, $I$ be a bounded closed interval of $\R$ and $C:t\in I \to C(t)$ be a map defined on $I$ taking values in the set of closed $\eta$-prox-regular subsets of $H$. Let us assume that $C(\cdot)$ varies in an absolutely continuous way, that is to say, there exists an absolutely continuous function $w:I\rightarrow \R$ such that, for any $y\in H$ and $s,t\in I$
\be{A1} \left|d(y,C(t))-d(y,C(s))\right|\leq |w(t)-w(s)|. \tag{A1} \ee
Let $f:I\times H \rightarrow H$ be a mapping which is measurable with respect to the first variable and continuous with respect to the second one and such that there exists a nonnegative function $\beta\in L^1(I,\R)$ satisfying for all $t\in I$ and for all $x\in \cup_{s\in I} C(s)$,
\be{A2} \|f(t,x)\|\leq \beta(t)\left(1+\|x\|\right). \tag{A2} \ee
Moreover, we suppose that $f$ satisfies a Lipschitz condition: for every $M>0$ there exists a non-negative function $k_M(.) \in L^1(I,\R)$ such that for all $t\in I$ and for any $(x,y)\in B(0,M)\times B(0,M)$,
\be{A3} \|f(t,x)-f(t,y)\|\leq k_M(t)\|x-y\|. \tag{L} \ee
Then for all $u_0\in C(0)$, the differential inclusion (\ref{eq:syst1}) has one and only one absolutely continuous solution.
\end{thm}

\begin{rem} \label{rem:mult} In \cite{Thibbv}, the authors describe more general results of existence for sweeping process with a multivalued perturbation: they deal with a perturbation $F:I\times H \rightrightarrows H$, which is assumed to be separately scalarly upper semicontinuous, admitting a compact and linear growth, and such that for all $x\in H$ the function $F(\cdot,x)$ has a measurable selection. We do not detail these assumptions as we will only consider the case of a single-valued mapping $f$.
\end{rem}

\mb We want to use the ``hypomonotonicity'' property of the proximal normal cone $\NN(C(t),.)$ to obtain information about the differential inclusion (\ref{eq:syst1}). First we describe a result concerning a constant set $C$.

\begin{prop} \label{prop1} Let $H$ be a Hilbert space,$C$ be a uniformly prox-regular subset and $f:I\times H \to H$ be a mapping satisfying the assumptions of Theorem \ref{thm1}. Then for all $u_0\in C$, the (unique) solution $u$ of (\ref{eq:syst1}) satisfies the following differential equation: for almost all $t_0\in I$, $u$ can be differenciate at the point $t_0$ and
\begin{equation} \label{eq:syst2}
\frac{du}{dt}(t_0) + P_{\NN(C,u(t_0))}\left[f(t_0,u(t_0))\right] = f(t_0,u(t_0)).
\end{equation}
\end{prop}

\dem For convenience and to expose the main arguments, we assume that $f$ is bounded on $I\times H$. In fact due to Assumption (\ref{A2}), we know that this property holds locally on time almost everywhere on $I$. As we are looking for local results, this restriction is allowed. \\
We follow the ideas of H. Brezis (see \cite{Brezis}), who has already proved similar results, in considering
multivalued maximal monotone operators instead of the proximal normal cone $\NN(C,u(t))$. The proximal normal cone
$\NN(C,\cdot)$ is not monotone, fortunately it is hypomonotone (a little weaker property) because of the uniform prox-regularity of the set $C$. By the work of R.A. Poliquin, R.T. Rockafellar and L. Thibault (see
\cite{PRT}), $\NN(C,.)$ satisfies: for all $z_1,z_2\in C$, $\zeta_1\in \NN(C,z_1)$ and $\zeta_2\in \NN(C,z_2)$ \be{eq:hypo} \langle \zeta_1-\zeta_2 , z_1-z_2
\rangle \geq - \frac{\|\zeta_1\|+\|\zeta_2\|}{\eta} \left\|z_1-z_2\right\|^2.  \ee 
From this property, we can obtain the desired result. \\
The function $\frac{du}{dt}$ belongs to
$L^{1}([0,T],H)$ so almost every point is a Lebesgue point of $\dot{u}$. 
The same reasoning holds for $t\to f(t,u(t))$, which is a bounded function.
Let $t_0\in I$ be a Lebesgue point for $\dot{u}$ and $ f(\cdot,u(\cdot))$. \\
Let us consider the following mapping $g$, defined on $H$ by
$$ g(v):=P_{\NN(C,v)}\left[f(t_0,v)\right].$$
For every point $t\in I$ and $v\in C$, the projection onto $\NN(C,v)$, due to its convexity, is everywhere well-defined and so
$P_{\NN(C,v)}\left[f(t_0,v)\right]$ corresponds to a unique point. \\
The constant function $\tilde{u}(t):=u(t_0)$ satisfies the following differential inclusion
\begin{equation} \frac{d\tilde{u}}{dt} + \NN(C,\tilde{u}) \ni g(\tilde{u}). \label{eqg} \ee
Let us first check that for all $t_0 <t$, we have:
\begin{equation} \label{eq:amont}
 \left\|u(t) - \tilde{u}(t) \right\| \leq \left\|u(t_0) - \tilde{u}(t_0) \right\| + \int_{t_0}^t \left[\left\|f(\sigma,u(\sigma))-g(\sigma,\tilde{u}(\sigma))\right\|+h(\sigma) \left\|u(\sigma) - \tilde{u}(\sigma) \right\|\right] d\sigma,
\end{equation}
where $h$ is given by
$$ h(\sigma):= \frac{1}{\eta} \left(\left\|\frac{du}{dt}(\sigma)-f(\sigma,u(\sigma))\right\| + \left\|\frac{d\tilde{u}}{dt}(\sigma)-g(\tilde{u}(\sigma))\right\| \right)\in L^1_{loc}(I).$$
Using both differential inclusions ((\ref{eq:syst1}) for $u$ and (\ref{eqg}) for $\tilde{u}$) and the
hypomonotonicity property of the proximal normal cone (\ref{eq:hypo}), we get:
\begin{align}
 \frac{1}{2} \frac{d}{dt} \left\|u(t) - \tilde{u}(t) \right\|^2 & = \left\langle \frac{du}{dt}(t) - \frac{d\tilde{u}}{dt}(t) , u(t)- \tilde{u}(t) \right\rangle \nonumber \\
 & \leq \left\langle f(t,u(t))-g(\tilde{u}(t)), u(t)- \tilde{u}(t) \right\rangle + h(t) \left\|u(t) - \tilde{u}(t) \right\|^2. \label{eq:hypo2}
\end{align}
The integration of this inequality on $[s,t]\subset I$ yields
\begin{align*}
\lefteqn{ \frac{1}{2}\left\|u(t) - \tilde{u}(t) \right\|^2 - \frac{1}{2}\left\|u(s) - \tilde{u}(s) \right\|^2} & & \\
 & &  \leq \int_{s}^t
\left[\left\|f(\sigma,u(\sigma))-g(\tilde{u}(\sigma))\right\|+h(\sigma) \left\|u(\sigma) - \tilde{u}(\sigma) \right\|\right] \left\|u(\sigma)- \tilde{u}(\sigma)\right\| d\sigma.
\end{align*}
Then we deduce (\ref{eq:amont}) with the help of Lemma A.5 in \cite{Brezis}. \\
Now we use that $\tilde{u}$ is constant and equal to $u(t_0)$. For $t=t_0+\epsilon$, we obtain
\begin{align*}
\lefteqn{\left\|u(t_0+\epsilon) - u(t_0) \right\|} & & \\
 & &  \leq \int_{t_0}^{t_0+\epsilon}
\left[\left\|f(\sigma,u(\sigma))-P_{\NN(C,u(t_0))}\left[f(t_0,u(t_0))\right]\right\|+h(\sigma)
\left\|u(\sigma) - u(t_0) \right\|\right] d\sigma. \end{align*} 
Finally we have
\begin{align*}
\limsup_{\epsilon\to 0} \frac{\left\|u(t_0+\epsilon) - u(t_0) \right\|}{\epsilon} & \leq \limsup_{\epsilon\to 0} \frac{1}{\epsilon} \int_{t_0}^{t_0+\epsilon} \left\|f(\sigma,u(\sigma))-P_{\NN(C,u(t_0))}\left[f(t_0,u(t_0))\right]\right\| d\sigma \\
 &  + \limsup_{\epsilon\to 0} \frac{1}{\epsilon} \int_{t_0}^{t_0+\epsilon} h(\sigma) \left\|u(\sigma) - u(t_0) \right\| d\sigma.
\end{align*}
The second term of the right member is vanishing as $t_0$ is a Lebesgue point of $h$ and $u$ is continuous at $t_0$. For the first term, we use that $t_0$ is a Lebesgue point of $f(\cdot,u(\cdot))$. It comes 
\begin{align*}
\lefteqn{\limsup_{\epsilon\to 0}
\frac{1}{\epsilon} \int_{t_0}^{t_0+\epsilon}
\left\|f(\sigma,u(\sigma))-P_{\NN(C,u(t_0))}\left[f(t_0,u(t_0))\right]\right\| d\sigma } & & \nonumber \\
 & & = \left\|f(t_0,u(t_0))-P_{\NN(C,u(t_0))}\left[f(t_0,u(t_0))\right]\right\|, 
\end{align*} 
and consequently
\begin{align} 
\limsup_{\epsilon\to 0} \frac{\left\|u(t_0+\epsilon) - u(t_0) \right\|}{\epsilon} & \leq \left\|f(t_0,u(t_0))-P_{\NN(C,u(t_0))}\left[f(t_0,u(t_0))\right]\right\| \nonumber  \\
 & \leq d(f(t_0,u(t_0)),\NN(C,u(t_0))) .  \label{eq:res}
 \end{align}
However we know that when $u$ is differentiable at $t$, then
$$\frac{du(t)}{dt} \in f(t,u(t))-\NN(C,u(t)).$$
Equation (\ref{eq:res}) gives us the desired
equality:
$$ \frac{du}{dt}(t_0) + P_{\NN(C,u(t_0))}\left[f(t_0,u(t_0))\right] = f(t_0,u(t_0)).$$
 \findem

\mb In this particular case of a constant prox-regular set $C$, we have obtained the equivalence between the differential inclusion (\ref{eq:syst1}) and the differential equation (\ref{eq:syst2}). In a general situation with a moving set $C(t)$, such equivalence may not hold (it is easy to build a counter-example with no perturbation $f=0$). 

\mb Using similar reasonning, we can describe a stability for the solutions of (\ref{eq:syst1}), already proved in Proposition 2 of \cite{Thibrelax}. We recall its proof for an easy reference.

\begin{prop} \label{cor} Under the assumptions of Theorem \ref{thm1}, for all $t\in I$ and $M$, there exists a constant $a>0$ (depending on $t$ and $M$) such that for the
solution $u$ (respectively $v$) associated to initial data $u_0$ (resp. $v_0$) with $\|u_0\|\leq M$ and $\|v_0\|\leq M$ we have:
$$ \left|u(t)-v(t) \right| \leq a \|u_0-v_0\|.$$
\end{prop}

\dem Let $u_0,v_0$ be two fixed data. Consider $u$ and $v$ a solution of (\ref{eq:syst1}) according to initial data
$u_0$ and $v_0$. Let $M'$ be a bound of the solution $u(t)$ and $v(t)$, depending only on $M$.
By the same reasoning (as for (\ref{eq:hypo2})) using the hypomonotonicity of the proximal normal
cone, we get:
$$ \frac{1}{2} \frac{d}{dt} \left\|u(t) - v(t) \right\|^2 \leq \left[k_{M'}(t)+ h(t) \right] \left\|u(t) - v(t) \right\|^2,
$$
where $h$ is defined as
$$ h(t):= \frac{1}{\eta} \left(\left\|\frac{du}{dt}(t)-f(t,u(t))\right\| + \left\|\frac{dv}{dt}(t)-f(t,v(t))\right\| \right)\in L^1(I).$$
Applying Gronwall's Lemma, we get
$$ \left\|u(t) - v(t) \right\| \leq \left\|u_0-v_0\right\| e^{\int_0^t \left[k_{M'}(\sigma)+h(\sigma)\right] d\sigma}.$$
Theorem 5.1 of \cite{Thibbv}) shows that we can estimate the function $h$ as folllows:
$$h(t) \leq \frac{2}{\eta}\left[ (1+M')\beta(t) + |\dot{w}(t)| \right]  \in L^1(I,\R).$$
As $k_{M'}\in L^1(I,\R)$, we also deduce the result. \findem

\mb Proposition \ref{prop1} gives an interesting result: for a non-moving set $C$, the important quantity seems to be $P_{\NN(C,u(t))}\left[f(t,u(t))\right]$, which is a particular point of the set $\NN(C,u(t))$. So we guess that we have not to require information for the whole cone $\NN(C,u(t))$ (obtained by the assumption of the uniform prox-regularity), but only on this specific point. This observation is the starting-point for the definition of ``directional prox-regularity'' (see Subsection \ref{subsec:dir}). 

\section{A particular example for a lack of ``uniform prox-regularity''. } \label{sec:ex}

The aim of this section is to describe with an example (due to the modelling of crowd motion in emergency evacuation) that the uniform ``prox-regularity'' of an interesting set may not be satisfied or it may be difficult to check this property. We refer the reader to the thesis of the second author (\cite{TheseJu}) for a complete and detailed description of this model and to \cite{MV,MV2} for articles.

\mb We quickly recall the model. It handles contacts, in order to deal with local interactions between people and to describe the whole dynamics of the pedestrian traffic. This microscopic model for crowd motion rests on two principles. On the one hand, each individual has a spontaneous velocity that he would like to have in the absence of other people. On the other hand, the actual velocity must take into account congestion. Those two principles lead to define the actual velocity field as the projection of the spontaneous velocity over the set of admissible velocities (regarding the non-overlapping constraints).

\mb We consider $N$ persons identified to rigid disks. For convenience, the disks are supposed here to have the same radius $r$. The
center of the $i$-th disk is denoted by $\qq_i\in \R^2$. Since overlapping is forbidden, the  
vector of positions $\qqq=(\qq_1,..,\qq_N) \in
\R^{2N}$ has to belong to the ``set of feasible configurations'', defined by 
\be{eq:Q} Q:=\left\{ \qqq  \in
\R^{2N},\ D_{ij}( \qqq) \geq 0\quad   \forall \,i \neq j  \right\},  \ee 
where $ D_{ij}(\qqq)=|\qq_i-\qq_j|-2r $ is the signed distance between disks $i$ and $j$. 

\mb We denote by $\UU(\qqq)=(U_1(\qq_1),..,U_N(\qq_N)) \in \R^{2N}$ the global spontaneous velocity of the crowd.
To get the actual velocity, we introduce the ``set of feasible velocities'' defined by:
 $$ \CCC_{\qqq}=\left\{  \vv  \in
\R^{2N} , \ \forall i<j \hspace{5mm} D_{ij}(\qqq)=0 \hspace{3mm}
\Rightarrow \hspace{3mm}
\GG_{ij}(\qqq)\cdot \vv \geq 0 \right\},  $$ with
$$\GG_{ij}(\qqq)=\nabla D_{ij}(\qqq)=
(0,\dots,0, -\eee_{ij}(\qqq),0,\dots,0,\eee_{ij}(\qqq),0,\dots,0 )\in
\R^{2N}  $$ and $\eee_{ij}(\qqq)= \frac{\qq_j-\qq_i}{|\qq_j-\qq_i|}$.
The actual velocity field is defined as the feasible field which is the closest to $\UU$ in the least square sense, which writes
\be{eq:model1} \frac{d\qqq}{dt} = P_{\CCC_{\qqq}}\left[\UU(\qqq)\right], \ee
where $P_{\CCC_{\qqq}}$ denotes the Euclidean projection onto the closed convex cone $\CCC_{\qqq}$.
Then using the Hilbertian structure of $\R^{2N}$ and convex analysis, we have the following results:

\begin{prop} \label{prop:model1}  The negative polar cone $\NNN_\qqq$ of $\CCC_\qqq $, i.e.,
$$ \NNN_\qqq := \CCC_\qqq^\circ
:=  \left\{\ww \in \R^{2N},\ \ww \cdot \vv \leq 0 \quad \forall \vv \in \CCC_\qqq \right\},$$
is equal to the proximal normal cone $N(Q,\qqq)$ and
$$ \NNN_\qqq = N(Q,\qqq) = \left\{-\sum \lambda_{ij} \GG_{ij}(\qqq),\ \lambda_{ij} \geq 0,\ 
D_{ij}(\qqq) > 0 \Longrightarrow \lambda_{ij} = 0
\right\} .$$
\end{prop}

\mb Using the classical orthogonal decomposition with two mutually polar cones (see \cite{JJM}), the main equation (\ref{eq:model1}) becomes
\be{eq:model2} \frac{d\qqq}{dt} + P_{N(Q,\qqq)}\left[\UU(\qqq)\right] = \UU(\qqq).\ee
As proved in Section \ref{sec1}, we know that for a Lipschitz bounded map $\UU$, this differential equation is equivalent to the following sweeping process:
\be{eq:model3} \frac{d\qqq}{dt} + N(Q,\qqq) \ni \UU(\qqq).\ee

\mb  The uniform prox-regularity of the set $C$ guarantees the existence and the uniqueness of solution for such a differential inclusion. This main property is proved in \cite{TheseJu,MV} in the case of this model where no obstacles have been taken into account:

\begin{thm} The set $Q\subset \R^{2N}$, defined by (\ref{eq:Q}) (which corresponds to a model without obstacles) is $\eta$-prox-regular with a constant $\eta=\eta(N,r)>0$.
 \end{thm}

\mb  We emphasize that this property was already quite difficult (to be proven) and use specific geometric properties, mainly precise estimates about the angles between the different vectors $\GG_{ij}(\qqq)$. \\
Now we are interested in extending this result with a model taking into account obtacles. We do not write details and hope to deal in a more precise way with this particular problem in a forthcoming work. We just want to explain with the following example how the assumption of ``uniform prox-regularity'' could fail.

\mb We consider a small parameter $\epsilon>0$ and two additional
obstacles represented by the lines $x=0$ and $x=4r-2\epsilon$ in the
physical plane and we consider two disks ($N=2$). The set of
feasible configurations $Q$ is now defined by
\be{eq:Qobs}
 Q:=\left\{ \qqq=(\qq_{1x},\qq_{1y},\qq_{2x},\qq_{2y}) \in \R^{4},\ D_{12}( \qqq) \geq 0\quad,\quad  r\leq \qq_{1x},\ \qq_{2x}\leq 3r-2\epsilon \right\}.  \ee 
We claim that if the set $Q$ is uniformly prox-regular then its constant has to be lower than $\sqrt{\epsilon}$. Indeed we can consider the specific configuration (represented in Figure \ref{figure1}) 
$$\qqq_0=(r-\epsilon,0,3r-\epsilon,0).$$

\begin{figure}[ht!]
\centering
\psfrag{a}[c]{$x=0$}
\psfrag{b}[c]{$r-\epsilon$}
\psfrag{c}[c]{$3r-\epsilon$}
\psfrag{d}[c]{$x=4r-2\epsilon$}
\includegraphics[width=0.4\textwidth]{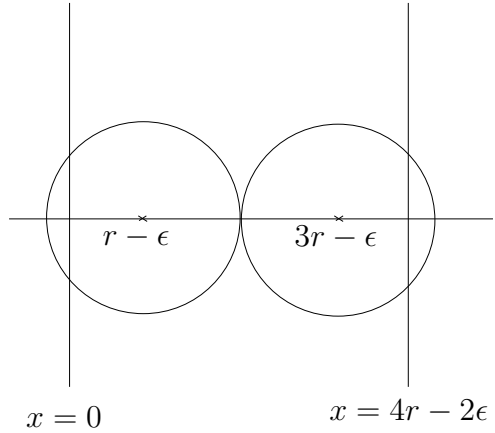}
\caption{Specific configuration}
\label{figure1}
\end{figure}

\mb The point $\qqq_0$ does not belong to $Q$ however $d(\qqq_0,Q)\lesssim \sqrt{\epsilon}$ and in invoking symmetry, it is obvious that this configuration does not admit a unique projection on $Q$. So if $Q$ is uniformly prox-regular then its constant of prox-regularity must be lower than $\sqrt{\epsilon}$. Furthermore similar configurations seem to produce some difficulties also in the numerical resolution. Indeed the Kuhn-Tucker multipliers (appearing in the discretization of the differential inclusion) could not be bounded in considering obstacles. We refer the reader to Remark 4.23 of \cite{TheseJu}.

\mb In conclusion, when we consider obstacles in the model of crowd
motion, the eventual uniform prox-regularity of $Q$ will depend on
the geometry of the obstacles, more precisely on their relative positions. This dependence is probably very
difficult to be estimated. Fortunately, as we are going to explain, it is not necessary to study the prox-regularity for all directions. Based on the proof of the existence of solutions and as we explain in the next sections, we only have to measure the prox-regularity in the direction given by $\UU(\qqq)$.

\mb Let us treat a special choice of $\UU$. For $\qq\in\R^2$, we define $U(\qq)$ as the unit vector directed by the shortest path avoiding obstacles from the point $\qq$ to the nearest exit (of the considered room) and then define 
$$\UU(\qqq)=(U(\qq_1),..,U(\qq_N)).$$
In Figure \ref{flotgeodesique} we consider a room containing obstacles with an exit (represented by the bold segment on the left). We draw the level curves of the distance function to the exit (obtained by a Fast Marching Method, see \cite{TheseJu}) and we represent the velocity field (corresponding to the gradient of this geodesic distance).

\begin{figure}[ht!]
\centering
 \includegraphics[width=0.46\textwidth]{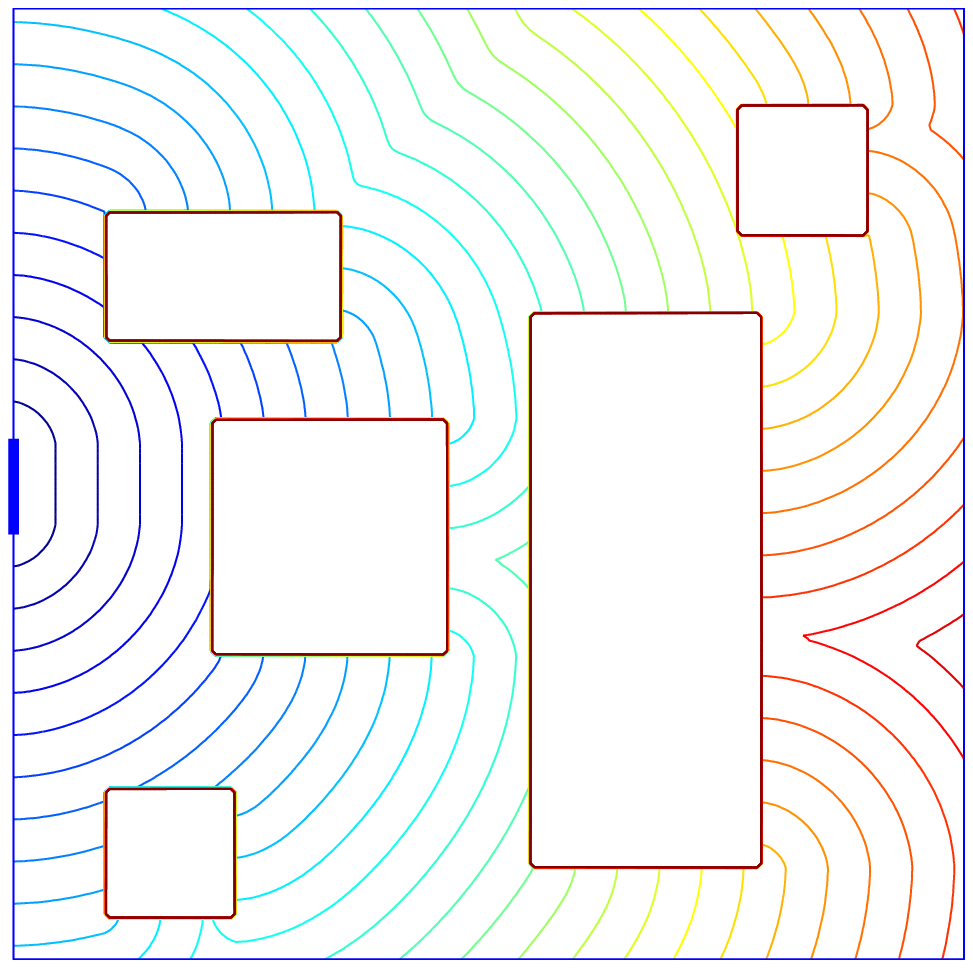} 
 \includegraphics[width=0.475\textwidth]{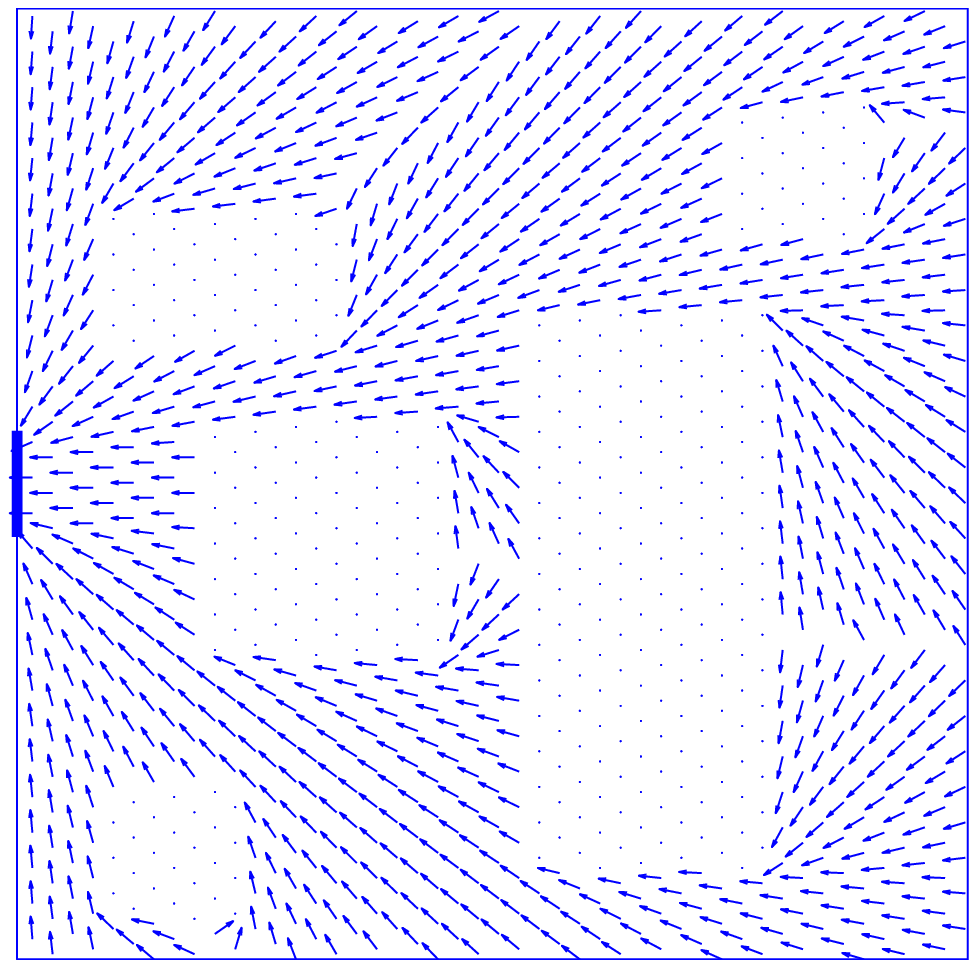} 
\caption{Level curves of the distance function and vector field of the velocity.}
\label{flotgeodesique}
\end{figure}

\mb We can also see that a person moving with velocity $\UU$ avoids the different obstacles. Indeed, with an elementary (infinitesimal) displacement in the direction given by this velocity field, we guess that the persons will not interact with the obstacles, as the velocity field get smoothly around the obstacles. The information about obstacles are now hidden in the vector field $\UU$. So we foresee that in the direction $\UU$, the new configuration $\qqq+h\UU(\qqq)$ (with a small enough parameter $h$) has only overlapping between disks (representing the persons) and no overlapping with obstacles. Consequently the configuration (illustrated by Figure \ref{figure1}) will never be realized by the crowd motion. 

\mb We do not give more details for this example, as it is not the aim of this paper. We would just like to emphasize that in this application, the uniform prox-regularity will not be easily checked but a kind of ``directional prox-regularity'' (along the perturbation $\UU(\qqq)$) would be more easily estimated. This is why, we propose in the next Section a rigourous definition of ``directional prox-regularity'' (motivated by this example) and then study sweeping process under this new assumption.

\section{About our assumptions ...}

We devote this section to the definitions of some new concepts needed in our assumptions. We want first to weaken the uniform ``prox-regularity'' assumption about the set $C$, in only requiring a ``directional prox-regularity''. Then we define a new property for the Banach space $\B$ (which generalizes a property of the Hilbertian structure) permitting us to prove existence for sweeping process.

\subsection{Concept of ``directional prox-regularity''.} \label{subsec:dir}

\mb Due to Proposition \ref{prop1}, we guess that it is not necessary to require the whole (in all the directions) property of uniform ``prox-regularity'' for the set $C(t)$. Indeed, during the construction of the solution (see the proof of Theorem \ref{thm1} and our proof of Theorem \ref{thmm}), we understand that we have to use the term
$$ P_{C(t)}\left[u(t)-hf(t,u(t))\right]$$
for  a small enough parameter $h$. This term is obviously well-defined for a uniformly prox-regular set $C(t)$.

\mb According to this observation, we define a new assumption, which only corresponds to a ``directional prox-regularity'', which will be sufficient to obtain an existence result of the differential inclusion.

\mb For more convenience, we will deal only with a simple case to introduce our concepts: we suppose that the set-valued map $C(\cdot)$ is constant. The case of a non constant set $C$ seems to be more technical as we need to know how the set $C$ is moving to only require a directional prox-regularity (see the comments after Theorem \ref{thmm2} and Theorem \ref{thmm3}).\\
Let $C$ be a fixed closed subset of a Banach space $\B$.

\begin{df} For every point $x\in C$ and $r>0$, we define $\Gamma^{r}(C,x)$ as the set of ``good directions $v$ to project at the scale $r$'' from $x+rv$ to $x$:
$$ \Gamma^{r}(C,x):= \left\{ v\in \B,\ x\in P_{C}(x+rv) \right\}.$$
\end{df}

\begin{rem} For all $x\in C$, we obviously have by definition of the proximal normal cone $\NN(C,x)=\cup_{r>0} \Gamma^r(C,x)$. 
\end{rem}

\begin{df} \label{df:Fprox} Let $f:\B \to \B$ be a mapping. We say that the set $C$ is ``$(r,F)$ prox-regular'' or ``$r$-prox-regular in the direction $f$'' if for all $x\in C$ and $s\in(0,r)$ 
\begin{itemize}
 \item[a)] the following projection is well-defined~:
 $$ z:=P_C\left(x+s\frac{f(x)}{\|f(x)\|}\right) $$
 \item[b)] and it satisfies
 $$ \frac{x+s\frac{f(x)}{\|f(x)\|}-z}{\left\| x+s\frac{f(x)}{\|f(x)\|}-z \right\|} \in \Gamma^r(C,z). $$
\end{itemize}
If $v=0$, we set $\frac{v}{\|v\|}:=0$ by convention.
\end{df}

\begin{rem} If the set $C$ is $r$-prox-regular then for all mappings $f$, it is $(r,f)$-prox-regular.
\end{rem}

\mb We can describe this definition as follows with Figure \ref{fig:fig3}.

\begin{figure}[ht!]
\centering
\psfrag{C}[c]{$C$}
\psfrag{x}[c]{$x \ \ $}
\psfrag{z}[c]{$z$}
\psfrag{y}[c]{$\ \ \ y$}
\psfrag{w}[c]{$\ \ \ w$}
\psfrag{F}[c]{$f(x) \hspace{0.3cm} $}
\includegraphics[width=0.4\textwidth]{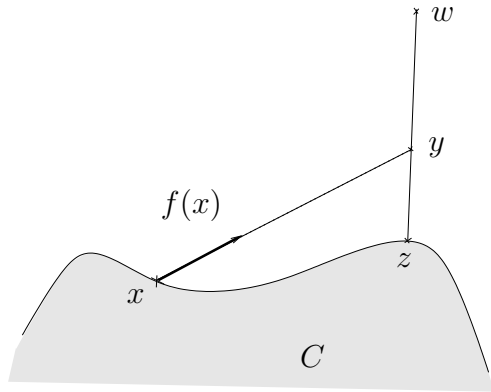}
\label{fig:fig3}
\caption{Illustration of the ``directional prox-regularity''.}
\end{figure}

\mb Let $x\in C$ be a point inside or on the boundary of set $C$ and let $y=x+s\frac{f(x)}{\|f(x)\|}$ (for $s\in(0,r)$) corresponding to a small perturbation of $x$ in the direction $f(x)$. We do not know if the point $y$ belongs to the set $C$ or not but we require that it stays in a {\em good} neighbourhood of set $C$. Referring to Condition $a)$, we ask that the projection of $y$ onto $C$ is well-defined, $z:=P_C(y)$. Consequently all points belonging to the segment $[y,z]$ project themselves on $z$. We have to be careful, as $r$ (equals to the length of $[x,y]$) could be larger than $d(y,C)$ we require with Condition $b)$ that all points belonging to segment $[w,z]$ (the segment of length $r$ extending the previous one) satisfy the same property. 

\mb We refer the reader to the work \cite{PRT} of  R.A. Poliquin, R.T. Rockafellar and L. Thibault. They define another concept of directional prox-regularity. Their notion is not comparable to our one as they only consider proximal normal directions. We just describe an interesting example~: in $\R^2$, let us consider the set $C$ defined by
$$ C:=\left\{(x,y)\in\R^2,\ x\leq 0 \textrm{ or } y\leq 0\right\}.$$
It is well-known that the point $0\in C$ is not regular and have no non-zero proximal directions. So at this point, the set $C$ is not prox-regular in any direction in the sense of \cite{PRT}. However the set $C$ is $r$-prox-regular in the direction $f$ for some mappings $f$. For example, it is easy to see that the set $C$ is $\infty$-prox-regular in the direction $f$ with $f(x,y):=(-1,-1)$ for all $(x,y)\in \R^2$. \\
This example shows how this new concept of directional prox-regularity can be far more weaker than the uniform prox-regularity.

\mb Now we have to define a new concept for Banach space which will ensure the existence of solutions for sweeping process. This is the goal of the next subsection.

\subsection{Geometry of Banach spaces} \label{subsec:geo}

First we recall some useful definitions, due to the geometric theory of Banach spaces (we refer the reader to
\cite{Diestel} for these concepts and more details).

\begin{df} Let $E$ be a Banach space, equipped with its norm $\|\ \|_{E}$.
\begin{itemize}
 \item The space $E$ is said to be {\em uniformly convex} if for all $\epsilon>0$, there is some $\delta>0$ so that for
any two vectors $x,y\in E$ with $\|x\|_E\leq 1$ and $\|y\|_E\leq 1$ we have
$$ \|x+y\|_E >2-\delta \Rightarrow \|x-y\|_E \leq \epsilon.$$
 \item The space $E$ is said to be {\em uniformly smooth} if the norm is uniformly {\em Fr\'echet differentiable} away
 of $0$, it means that for any two unit vectors $x_0,h\in E$, the limit
 $$ \lim_{t\to 0} \frac{\|x_0+th\|_E-\|x_0\|_E}{t}$$
 exists uniformly with respect to $h,x_0\in S(0,1)$.
\end{itemize}
We write $S(0,1):=\{x\in E,\ \|x\|_E=1\}$ for the unit sphere.
\end{df}

\mb We refer to \cite{BTZ} (Lemma 2.1) for the following geometric lemma:

\begin{lem} \label{lem:geo2} Let $\B$ be a Banach space and $C$ be a closed subset of $\B$. Then for $x\in C$ and $v\in \Gamma^r(x)$, we have $\lambda v\in\Gamma^r(x)$ for all $\lambda\in(0,1)$. Therefore if we assume that $\B$ is uniformly convex then for all $\lambda\in(0,1)$, we have $x=P_C(x+\lambda r v)$.
\end{lem}

\mb The first part is well-known (see for example Property 2.19 of \cite{TheseJu}), the second part is quite more complicated.

\mb We recall a famous result (due to D. Milman and B.J. Pettis \cite{M,P})

\begin{thm} If a Banach space $E$ is uniformly convex then it is reflexive: $E^{**}=E$.
\end{thm}

\mb The following well-known results (see e.g., \cite{Diestel}) will be also needed:

\begin{thm} If $E^*$ is uniformly convex then $E$ is uniformly smooth and $E$ is reflexive.
If $E^*$ is separable then $E$ is separable.
\end{thm}

\mb Now we consider some results concerning the smoothness of the norm.

\begin{rem} If $E$ is uniformly smooth, then $x\rightarrow \|x\|_E$ is $C^1$ on $E\setminus \{0\}$.
\end{rem}

\begin{prop} \label{prop:norm1} If $E$ is uniformly smooth, then for all $x\in E\setminus \{0\}$, we have
$$ \langle(\nabla\|.\|_E)(x),x\rangle = \|x\|_E.$$
By triangle inequality, $\left\|(\nabla\|.\|_E)(x)\right\|_{E^*}=1$.
\end{prop}

\mb Now as we know that the norm could be non differentiable at the origin $0$, we study the function $x\to \|x\|_E^p$
for an exponent $p>1$.

\begin{prop} \label{prop:JC1} Let $E$ be a uniformly smooth Banach space and $p\in(1,\infty)$ be an exponent. The function $x\to \|x\|_E^p$ is $C^1$ over
the whole space $E$.
\end{prop}

\mb For an easy reference, we explain the proof: \\

\dem As the norm is $C^1$ on $E\setminus \{0\}$, we have just to check the claim at the point $0$. As for every $h\in
E$, $\frac{\|th\|_E^{p}}{t}$ tends to $0$ when $t\to 0$, we deduce that $\phi:=\|.\|_E^p$ is differentiable at $0$ and
its gradient is null at this point. We have now to check that $\nabla \phi$ is continuous at this point. For any non
zero vector $x$, using Proposition \ref{prop:norm1}, we get:
$$\left\| \nabla \phi(x) \right\|_{E^*} \leq p \|x\|_{E}^{p-1} \xrightarrow{x\to 0} 0.$$
So we have proved that $\nabla \phi$ is continuous at $0$, which concludes the proof. \findem

\begin{df} For $E$ an uniformly smooth Banach space and $p\in(1,\infty)$, we denote
$$J_p(x) := \frac{1}{p} \left( \nabla \|.\|_{E}^p \right)(x) \in E^*.$$
\end{df}

\begin{rem} These mappings were already appeared to study the prox-regularity of a set (for example) in the work of F.~Bernard, L.~Thibault and N.~Zlateva (see \cite{BTZ,BTZ2}). We refer the reader to the work \cite{XR} of Z.B. Xu and G.F. Roach for more details about these mappings in an abstract framework.
\end{rem}

\begin{prop} \label{prop:unif} Let $E$ be a uniformly smooth Banach space and $p\in[2,\infty)$ be an exponent. Then $J_p$ is locally uniformly continuous: for all $\epsilon>0$, there exists $\delta>0$ such that for all $x,y\in E$
\be{ee} \left. \begin{array}{l}
\|x\|_E \leq 1 \vsp \\
\|y\|_{E}\leq 1 \vsp \\
 \|x-y\|_{E} \leq \delta
\end{array} \right\} \Longrightarrow \left\|J_p(x)-J_p(y) \right\|_{E^*} \leq \epsilon. \ee
 \end{prop}

\dem Just for convenience, we deal only with $p=2$. \\
Since the space is uniformly smooth, we know that $J_1(x)$ is uniformly continuous near $S(0,1)$ (see \cite{Diestel}). So let $\epsilon$ be fixed, we recall $J_2(x):=\|x\|_E J_1(x)$ and take $\eta$ such that if $\|z-z'\|_E\leq \eta$ then $\|J_1(z)-J_1(z')\|_{E^*}\leq \epsilon$ for $z,z'\in B(0,2)\setminus B(0,1/2)$. We set $\delta= min\{\epsilon/9,\eta \}$. \\
Take two points $x,y$ satisfying the assumption of (\ref{ee}). If $\|x\|_{E}\leq 4\delta$ then $\|y\|_{E}\leq 5\delta$ and so as $J_1$ is bounded 
$$ \left\|J_2(x)-J_2(y) \right\|_{E^*} \leq 9\delta+\epsilon\leq 2\epsilon.$$
Assume now that $\|x\|_{E}\geq 4\delta$ then with $\lambda=\|x\|_{E}^{-1}\geq 1$, we have $J_2(\lambda x)=\lambda J_2(x)$
so
$$ \left\|J_2(x) - J_2(y) \right\|_{E^*} \leq \frac{1}{\lambda}\left\| J_2(\lambda x) - J_2(\lambda y) \right\|_{E^*}\leq \|x\|_{E} \left\| J_2(\lambda x) - J_2(\lambda y) \right\|_{E^*}.$$
Now the whole segment $[\lambda x,\lambda y]$ is included in a neighbourhood of the sphere $S(0,1)$ and we have
$$ \left\|\lambda x - \lambda y \right\|_{E} \leq \lambda \|x-y\|_{E} \leq \lambda \delta.$$
We can also divide $[\lambda x,\lambda y]$ by $\lfloor\lambda\rfloor+1$ intervals of length $\delta$ (all of them included in a neighbourhood of the corona). Using the uniform continuity of $J_2$ around the sphere and $\lambda\geq 1$, we deduce
$$ \left\| J_2(\lambda x) - J_2(\lambda y) \right\|_{E^*} \leq (\lambda+1)\epsilon \leq 2 \frac{\epsilon}{\|x\|_{E}},$$
which permits us to deduce the desired inequality. \findem

\mb Now we can describe the useful assumption:

\begin{df} \label{df:imp} Let $I$ be an interval of $\R$. A separable reflexive uniformly smooth Banach space $E$ is said to be ``$I$-smoothly weakly compact'' for an exponent $p\in(1,\infty)$ if for all bounded sequence
$(x_n)_{n\geq 0}$ of $L^\infty(I,E)$, we can extract a subsequence $(y_n)_{n\geq 0}$ weakly converging to a point $y\in
L^\infty(I,E)$ such that for all $z\in L^\infty(I,E)$ and $\phi\in L^1(I,\R)$,
\begin{align}
\lefteqn{\lim_{n\to\infty} \int_I {}_{E^*}\langle J_p(z(t)+y_n(t))-J_p(y_n(t)) , y_n(t)\rangle_E \ \phi(t) dt =} & & \nonumber \\
 & & \int_I  {}_{E^*}\langle J_p(z(t)+y(t))-J_p(y(t)) , y(t)\rangle_E\ \phi(t) dt.\label{hyp}
\end{align}
\end{df}

\begin{rem} It is easy to check that the notion of ``$I$-smoothly weak
 compactness'' does not depend on the time interval $I$.
\end{rem}

\begin{rem} 
As $E$ is reflexive and separable, $L^{\infty}(I,E)=\left[L^1(I,E^*)\right]^{*}$ and by the Banach-Alaoglu-Bourbaki Theorem, we know that we can extract a weak converging subsequence $(y_n)_n$ from the initial bounded sequence $(x_n)_n$. However this weak convergence is not sufficient to insure (\ref{hyp}) in general. 
\end{rem}

\mb First we give several examples to illustrate this definition and to show that it has a ``non trivial'' sense.

\begin{prop} \label{prop:hilbert} All separable Hilbert space $H$ are $I$-smoothly weakly compact for $p=2$. 
\end{prop}

\dem It is well-known that for a Hilbert space, $J_2$ is given by $J_2(x)=x$. So (\ref{hyp}) corresponds to \be{ex1}
 \lim_{n\to\infty} \int_I \phi(t)\langle z(t) , y_n(t)\rangle dt = \int_I \phi(t) \langle z(t), y(t)\rangle dt.\ee As
$L^\infty(I,H)=\left[L^1(I,H)\right]^*$, we know that we can find a subsequence $(y_n)_n$ which weakly converges to a
point $y\in L^\infty(I,H)$. In considering $\phi(\cdot)z(\cdot)\in L^1(I,H)$, we conclude the proof. \findem

\mb We can not prove that the whole Lebesgue spaces or the whole Sobolev spaces are $I$-smoothly weakly compact for an exponent. However under an extra constraint over the sequence $(y_n)_n$, the desired conclusion holds:

\begin{prop} Let $U$ be an open subset of $\R^n$ or a Riemannian manifold. For all even integer $p\in[2,\infty)$ and $s\geq 0$, the Sobolev space $E=W^{s,p}(U,\R)$ is $I$-smoothly weakly compact for $p$ under an extra assumption :  from any bounded sequence $(x_n)_n$ of $L^\infty(I,E)$, which is also bounded in $L^\infty(I,W^{s+1,p}(U))$, then there exists a subsequence $(y_n)_{n}$ satisfying (\ref{hyp}).
\end{prop}

\dem Just for convenience we deal with $s=1$ (else we have to use properties of the singular operator $(1-\Delta)^{-s/2}$).
In this case, we consider  a bounded sequence $(x_n)_n$ of $L^\infty(I,W^{1,p}(U))$. We leave to the reader the computation of the gradient $J_p$ and we claim that for $f\in W^{1,p}(U)$
$$\left\langle J_p(f),h\right\rangle= \langle f^{p-1},h\rangle  +\sum_{i=1}^n \left\langle \left(\frac{\partial f}{\partial x_i}\right)^{p-1}, \frac{\partial h}{\partial x_i} \right\rangle.$$ So to check (\ref{hyp}), it suffices to prove that there exists a subsequence $(y_n)_n$ weakly converging to $y\in L^\infty(I,W^{1,p}(U))$ such that for all $g\in L^\infty(I,W^{1,p}(U))$ and $\phi\in L^1(I,\R)$,
\begin{align}
 \lefteqn{\lim_{n\to\infty} \int_{I\times U}
\left[(g(t,x)+y_n(t,x))^{p-1}-(y_n(t,x))^{p-1}\right] y_n(t,x) \phi(t) dtdx } & & \nonumber \\
 & & =\int_{I\times U}
\left[(g(t,x)+y(t,x))^{p-1}-(y(t,x))^{p-1}\right] y(t,x) \phi(t) dtdx \label{ex2} \end{align}
 and for $i\in\{1,..,n\}$
 \begin{align}
 \lefteqn{ \lim_{n\to\infty} \int_{I\times U} \left[(\partial_{x_i} g(t,x)+\partial_{x_i}
y_n(t,x))^{p-1} -(\partial_{x_i} y_n(t,x))^{p-1}\right] \partial_{x_i} y_n(t,x) \phi(t) dtdx } & & \nonumber \\
& & =\int_{I\times U} \left[(\partial_{x_i} g(t,x)+\partial_{x_i} y(t,x))^{p-1}-(\partial_{x_i} y(t,x))^{p-1}\right]
\partial_{x_i} y(t,x) \phi(t) dtdx.\label{ex3}
\end{align}
 As $p$ is an integer, using the ``binomial formula'', we get:
$$ y_n\left[(g+y_n)^{p-1}-(y_n)^{p-1}\right] = \sum_{k=0}^{p-2} \begin{pmatrix} p-1 \\ k \end{pmatrix} {y_n}^{k+1} g^{p-1-k}.$$
Which is interesting is that $g(t,.)\in L^{p'}$ implies $g(t,.)^{p-1-k}\in L^{(p/(k+1))'}$ and $y_n(t,.)\in L^p$ implies
$y_n(t,.)^{k+1} \in L^{p/(k+1)}$. So from the initial bounded sequence $(x_n)_{n\geq 0}$, we know that we can extract a
subsequence $(y_n)_n$ which weakly converges to a function $y\in L^\infty(I,W^{1,p})$ and such that for all
$k\in\{0,..,p-2\}$, $(y_n^{k+1})_n$ weakly converges to $y^{k+1}$ in $L^\infty(I,W^{1,p/(k+1)})$. According to a particularity of
Sobolev spaces, we know that we can extract of sequences $(y_n^{k+1})_n$ a subsequence which weakly and almost everywhere
converges to $y_{k+1}\in L^\infty(I,W^{1,p/(k+1)})$ as $W^{1,p/(k+1)}$ is a reflexive space. Then we deduce that almost everywhere
$y_{k+1}=y^{k+1}$, this well-known property of Sobolev spaces (weak convergence implies almost everywhere convergence) was already studied, see \cite{MT} for example. Then as we have a finite sum of limits, we get (\ref{ex2}). \\
Similarly as the sequence $(y_n)_n$ is assumed to be bounded in $L^\infty(I,W^{2,p}(U))$, we can produce the similar reasoning and prove the limit (\ref{ex3}), which concludes the proof. \findem

\begin{prop} For all even integer $p\in[2,\infty)$, the Lebesgue space $l^{p}(\Z)$ is $I$-smoothly weakly compact for $p$.
\end{prop}

\mb  We leave the proof to the reader, it is easier than the previous one. The important fact here, is that we are working
on $\Z$, a discrete space. So a weakly convergent sequence converges pointwise everywhere.

\section{Study of sweeping process in an abstract framework.} \label{sec:res}

Sweeping process have been studied in numerous papers in the case of the euclidean space first and then in a Hilbert space. The main technical difficulty is to obtain a kind of ``weak
continuity'' of the projection $P_C$. This problem is solved because the support function of Clarke's subdifferential of the distance function $d(.,C)$ is upper semicontinuous, when $C$ is a uniformly prox-regular set.

\mb We propose here new arguments to get around this difficulty. These ones permit us to understand the useful assumptions on the Banach space which are required to obtain a result of existence.

\mb The following proposition describes this useful property: a kind of ``weak continuity of the map $x\to \Gamma^r(C,x)$''. We recall that $I$ corresponds to the bounded time interval.

\begin{prop} \label{prop:imp} Let $(\B,\|\ \|)$ be a separable, reflexive and uniformly smooth Banach space. Let $C \subset \B$ be a closed subset. We assume that for an exponent $p\in[2,\infty)$ and  a bounded sequence $(v_n)_{n\geq 0}$ of $L^\infty(I,\B)$, we can extract a subsequence $(v_{k(n)})_{n\geq 0}$ weakly converging to a point $v\in
L^\infty(I,\B)$ such that for all $z\in L^\infty(I,\B)$ and $\phi\in L^1(I,\R)$,
\begin{align}
\lefteqn{\limsup_{n\to\infty} \int_I {}_{\B^*}\langle J_p(z(t)+v_{k(n)}(t))-J_p(v_{k(n)}(t)) , v_{k(n)}(t)\rangle_\B \phi(t) dt } & & \nonumber \\
 & & \leq \int_I  {}_{\B^*}\langle J_p(z(t)+v(t))-J_p(v(t)) , v(t)\rangle_\B \phi(t) dt.\label{prop:hyp}
\end{align}
Then the projection $P_C$ is weakly continuous in $L^\infty(I,\B)$ (relatively to the directions given by the sequence $(v_n)_n$) in the following sense: for all $r>0$ and for any  bounded sequence $(u_n)_n$ of $L^\infty(I,C)$ satisfying 
$$ \left\{ \begin{array}{l}
            u_n \longrightarrow u \quad \textrm{in}\ L^\infty(I,\B) \vsp \\
            u_n(t)\in P_C(u_n(t)+rv_n(t))\quad \textrm{a.e.}\ t\in I,
           \end{array} \right. $$
one has for almost every $t\in I$
$$ u(t)\in P_C(u(t)+rv(t)). $$
\end{prop}

\mb The above assumption is satisfied if the Banach space $\B$ is supposed to be ``$I$-smoothly weakly compact'' for an exponent $p\in[2,\infty)$. We can rewrite the conclusion as follows if for all $t\in I$, $v_n(t)\in \Gamma^r(C,u_n(t))$, then at the limit it holds that  $v(t)\in \Gamma^r(C,u(t))$, for almost every $t\in I$.

\begin{rem} We emphasize that this proposition has no link with the prox-regularity of the set $C$. This property is purely topological and only depends on the considered Banach space $\B$.
 \end{rem}

\dem  With the homogeneity of $J_p$ ($J_p(sx)=s^{p-1}J_p(x)$), in replacing $sz(t)$ by $z(t)$ in (\ref{prop:hyp}), we have for all $s\in(0,r)$
\begin{align}
\lefteqn{\limsup_{n\to\infty} \int_I {}_{\B^*}\langle J_p(z(t)+sv_{k(n)}(t))-J_p(sv_{k(n)}(t)) , v_{k(n)}(t) \rangle_\B \phi(t) dt } & & \nonumber \\
 & & \leq \int_I  {}_{\B^*}\langle J_p(z(t)+sv(t))-J_p(sv(t)) , v(t)\rangle_\B \phi(t) dt.\label{prop:hyp2}
\end{align}
It remains to prove that for almost every $t\in I$, $v(t)\in \Gamma^r(C,u(t))$.
Fixing any $\xi\in C$, for all integer $n$ and almost every $t\in I$, as $u_{k(n)}(t)\in P_C(u_{k(n)}(t)+rv_{k(n)}(t))$, we have
$$ \left\| u_{k(n)}(t)+rv_{k(n)}(t)-\xi \right\|^p-\left\|rv_{k(n)}(t)\right\|^p\geq 0.$$
Using Proposition \ref{prop:JC1}, this inequality can be written $$ \int_0^r \frac{d}{ds} \left[\left\| u_{k(n)}(t)+sv_{k(n)}(t)-\xi \right\|^p-\left\|s
v_{k(n)}(t)\right\|^p\right]ds \geq -\left\|u_{k(n)}(t)-\xi \right\|^p $$ and so $$\int_0^r \left\langle
J_p(u_{k(n)}(t)+s v_{k(n)}(t)-\xi)-J_p(sv_{k(n)}(t)),v_{k(n)}(t) \right\rangle ds \geq -\frac{1}{p}\left\|u_{k(n)}(t)-\xi \right\|^p.$$ 
Then for all smooth nonnegative function $\phi\in L^1(I,\R)$, we have 
\begin{align}
 \int_0^r \int_I
\phi(t)\left\langle J_p(u_{k(n)}(t)+sv_{k(n)}(t)-\xi)-J_p(sv_{k(n)}(t)),v_{k(n)}(t) \right\rangle dtds  \nonumber \hspace{2cm}  \nonumber \\
 \hspace{6cm} \geq -\frac{1}{p}\left(\int_I
\phi(t)\left\|u_{k(n)}(t)-\xi \right\|^p dt\right). \label{ineg1}
\end{align} 
We are now looking for passing to the limit in this inequality in order to get 
\begin{align}
 \int_0^r \int_I
\phi(t)\left\langle J_p(u(t)+s v(t)-\xi)-J_p(sv(t)),v(t) \right\rangle dtds \hspace{2cm} \nonumber \\
 \hspace{6cm} \geq -\frac{1}{p}\left(\int_I
\phi(t)\left\|u(t)-\xi \right\|^p
dt\right). \label{ineg2}
\end{align}
 As $(u_n)_n$ is bounded in $L^\infty(I,\B)$ and strongly converges to $u$ in $L^\infty(I,\B)$, it is obvious that
$$ \lim_{n\to \infty} \left(\int_I \phi(t)\left\|u_{k(n)}(t)-\xi \right\|^p
dt\right) = \left(\int_I \phi(t)\left\|u(t)-\xi \right\|^p dt\right).$$ Now consider the left-side of
(\ref{ineg1}). We know that $J_p$ is always locally bounded in $\B^*$ and is locally uniformly continuous as
$\B$ is uniformly smooth (see Proposition \ref{prop:unif}). For almost every $t\in I$ and all $s\in[0,r]$, we have
$$ \lim_{n\to\infty} \left\langle J_p(u_{k(n)}(t)+s
v_{k(n)}(t)-\xi)-J_p(u(t)+sv_{k(n)}(t)-\xi),v_{k(n)}(t) \right\rangle  =0$$ and this convergence is uniform with respect to $t\in I$ and $s\in(0,r)$. So the limit and the integrals can be inverted (according to Lebesgue's Theorem) and then
\begin{align*}
 \lim_{n\to \infty} \left|\int_0^r \int_I \phi(t)\left\langle J_p(u_{k(n)}(t)+s
v_{k(n)}(t)-\xi)-J_p(sv_{k(n)}(t)),v_{k(n)}(t) \right\rangle dtds \right. \hspace{1cm} \\
 \left. \hspace{1cm}  -\int_0^r \int_I \phi(t)\left\langle J_p(u(t)+s
v_{k(n)}(t)-\xi)-J_p(sv_{k(n)}(t)),v_{k(n)}(t) \right\rangle dtds\right| =0.
\end{align*}
From (\ref{prop:hyp2}) with $z(t)=u(t)-\xi$ and by Fatou's Lemma\footnote{Although the quantities are not necessary nonnegative, Fatou's Lemma can be applied. This is due to the fact that the integrated quantity is bounded by a constant (only depending on $\|\phi\|_{L^1}$, $\|\xi\|$ and the two bounded sequences), which is obviously integrable on $I\times [0,r]$.}, we obtain
\begin{align}
 \limsup_{n\to \infty} \int_0^r \int_I \phi(t)\left\langle J_p(u(t)+s
v_{k(n)}(t)-\xi)-J_p(sv_{k(n)}(t)),v_{k(n)}(t) \right\rangle dtds \hspace{1cm} \nonumber \\
 \leq \int_0^r \int_I \phi(t)\left\langle J_p(u(t)+sv(t)-\xi)-J_p(sv(t)),v(t) \right\rangle dtds. \label{ineg3}
\end{align}
With (\ref{ineg1}) and (\ref{ineg3}), we can conclude the proof of (\ref{ineg2}). \\
Now we produce the inverse reasoning in integrating the gradient $J_p$, obtaining from (\ref{ineg2}) that
\be{inegphi} \int_I \phi(t) \left[\left\| u(t)+rv(t)-\xi \right\|^p-\left\|rv(t)\right\|^p\right] dt \geq 0. \ee
That holds for every nonnegative smooth function $\phi\in L^1(I,\R)$, so we deduce that there exists a measurable set
$A_\xi\subset I$ satisfying  $|A_\xi|=0$ and such that for all $t\in I\setminus A_\xi$
$$ \left\| u(t)+rv(t)-\xi \right\| \geq \left\|rv(t)\right\|.$$
Now we use that $\B$ is separable and so $C$ is too. By taking  a dense sequence $(\xi_i)_{i\geq 0}$ of $C$,
we define $A:=\cup_{i\geq 0} A_{\xi_i}$. Then $|A|=0$ and for all $t\in I\setminus A$ and all $i\geq 0$, we have
$$ \left\| u(t)+rv(t)-\xi_i \right\| \geq \left\|rv(t)\right\|.$$
This last inequality is continuous with respect to $\xi_i$ and so by density, holds for all $\xi\in C$. That proves
$$ u(t) \in P_{C}(u(t)+rv(t))$$
and concludes the proof.
  \findem

\begin{rem} \label{rem:xi} For the proof, we have used a constant point $\xi\in C$. We emphasize that the different arguments hold with a bounded time-measurable map $\xi(\cdot)$ defined on $I$ taking values in $C$ and permit to obtain (\ref{inegphi}). Then we have to use the separability of the space $L^\infty(I,\B)$ for the $L^1(I,\B)$-norm in order to complete the proof.
\end{rem}

\mb Now we are going to use this preliminary and technical result to study existence and uniqueness of sweeping process. We first describe a result of existence in some $I$- smoothly weakly compact Banach spaces. Then we give a more precise result in a Hilbert space and obtain uniqueness of the solution.

\subsection{Sweeping process in Banach spaces for a constant set $C$.}

\mb In the case of general Banach spaces, an extra assumption about the set $C$ will be required. We introduce this one ~:

\begin{df} A subset $C\subset \B$ is said to be {\em ball-compact} if for all closed ball $\overline{B}=\overline{B}(x,R)$, the set $\overline{B} \cap C$ is compact.
\end{df}

\mb Obviously a ball-compact subset $C$ is closed.

\mb We now come to our main result in a Banach space $\B$.

\begin{thm} \label{thmm} Let $I=[0,T]$ be a bounded time-interval and $\B$ be a separable, reflexive, uniformly smooth Banach space, which is ``$I$-smoothly weakly compact'' for an exponent $p\in[2,\infty)$.
Let $f: \B \to \B$ be a bounded and continuous function, $r>0$ and $C\subset \B$ a nonempty, ball-compact and $(r,f)$-prox-regular set. Then for all $u_0\in C$, the system \be{eq:syst33} \left\{ \begin{array}{ll}
 \dot{u}(t)+\NN(C,u(t)) \ni f(u(t)) \vsp \\
 u(0)=u_0
\end{array} \right. \ee
has an absolutely continuous solution $u$, which lives in $C$. Moreover we have for almost every $t\in I$
\be{eq:cond} \left\| \dot{u}(t)-f(u(t)) \right\|_\infty \leq \|f\|_\infty. \ee
\end{thm}

\mb Indeed we are going to solve the following stronger system:
\be{eq:sys33} \left\{ \begin{array}{ll}
 \dot{u}(t)+\Gamma^{r/\|f\|_\infty}(C,u(t)) \ni f(u(t)) \vsp  \\
 u(0)=u_0 \ .
\end{array} \right. \ee

\mb The proof is a mixture of the classical one (see the papers cited in the introduction) based on the construction of discretized solutions and of Proposition \ref{prop:imp} which permits us to study the limit function.

\dem We follow the ideas of the well-known proof using an uniform prox-regular set $C$ (see \cite{Thibnonconv}). For an easy reference, we recall it and we will emphasize why our assumption is sufficient. 

\mb {\bf First step: } Construction of ``discretized solutions''. \\
We fix a small enough scale $h=T/n$ such that
\be{eq:1} h \|f\|_\infty \leq r/2.  \ee
Consider a partition of the time-interval $I=[0,T]$ defined by $t_i^n = ih$ for $i\in\{0,..,n\}$. We build $(u_n^i)_{0\leq i\leq n}$ as follows:
\be{aez} \left\{ \begin{array}{l}
                u_n^0=u_0 \vsp \\
                u_n^{i+1}=P_{C}\left[u_n^i +hf(u_n^i)\right].
           \end{array} \right. \ee
This operation is allowed as $\left\| hf(u_n^i) \right\|\leq r/2$ and the set $C$ is assumed to be $r$-prox-regular in the direction $f$ (see Definition \ref{df:Fprox}).
Now we use the points $(u_n^i)_{0\leq i\leq n}$ to obtain two piecewise maps $u_n$ and $f_n$ on $I$ (taking values in $\B$) in defining their restriction to each interval $I_i:=[ih, (i+1)h[$ by putting for every $t\in I_i$:
$$ f_n(t):=f(u_n^i)$$
and
\begin{align*}
u_n(t) &:= u_n^i + \frac{t-ih}{h}\left[u_n^{i+1}-u_n^i-hf(u_n^i) \right] + \left[t-ih \right]f(u_n^i) \\
 & = u_n^i + \left(\frac{t}{h}-i\right)\left[u_n^{i+1}-u_n^i\right].
\end{align*}
The function $u_n$ is continuous on $[0,T]$.

\mb {\bf Second step: } Differential inclusion for the ``discretized solution''. \\
We look for a differential inclusion satisfied by the function $u_n$. For almost every $t\in I_i$, we have
$$ \frac{d u_n(t)}{dt} = \frac{1}{h}\left[u_n^{i+1}-u_n^i-hf(u_n^i) \right] + f_n(t).$$
We define $\Delta_n(t)$ as follows
$$\Delta_n(t):=\frac{d u_n(t)}{dt} - f_n(t) = \frac{1}{h}\left[u_n^{i+1}-u_n^i-hf(u_n^i) \right].$$
We claim that $-\Delta_n(t) \in \Gamma^{r/\|f\|_\infty}(C,u_n^{i+1})\cap B(0,\|f\|_\infty)$, which is equivalent to
\be{eq:discrete} \|\Delta_n(t)\|\leq \|f\|_\infty \qquad \textrm{and} \qquad P_{C}\left[u_n^{i+1} -\frac{r}{\|f\|_\infty} \Delta_n(t)\right] \ni u_n^{i+1}. \ee
First we check that $\Delta_n(t)$ is a bounded vector. Using the construction of the point $u_n^{i+1}$ and the fact that $u_n^i\in C$, we have
\begin{align}
 \left\|\Delta_n(t) \right\| & = \frac{1}{h}\left\| P_{C}\left[u_n^i +hf(u_n^i)\right] - \left[u_n^i+hf(u_n^i) \right] \right\| \nonumber \\
 & \leq \frac{1}{h}\left\| u_n^i - \left[u_n^i+hf(u_n^i) \right] \right\| \\
 & \leq \left\| f(u_n^i) \right\| \leq \|f\|_\infty. \label{eq:2}
\end{align}
Then considering the vector $v:=u_n^i+hf(u_n^i)$, we have
\begin{align*}
 u_n^{i+1} -\frac{r}{\|f\|_\infty}\Delta_n(t) & = u_{n}^{i+1}-\frac{r}{h\|f\|_\infty }\left[u_n^{i+1}-u_n^i-hf(u_n^i) \right] \\
 & = P_C(v) - \frac{r}{h\|f\|_\infty}\left[ P_C(v)-v \right].
\end{align*}
Since $C$ is $r$-prox-regular in the direction $f$, we know that
$$ P_C\left(P_C(v) - \frac{r}{\|P_C(v)-v \|}\left[ P_C(v)-v \right]\right) \ni P_C(v)=u_n^{i+1}.$$ 
From (\ref{eq:2}), we deduce that $\|P_C(v)-v \| \leq h\|f\|_\infty$ and so with  the geometric Lemma \ref{lem:geo2}, we get
$$ P_C\left(P_C(v) - \frac{r}{h \|f\|_\infty}\left[ P_C(v)-v \right]\right) \ni P_C(v)=u_n^{i+1},$$
which concludes the proof of (\ref{eq:discrete}). For the discretized solution $u_n$, we have proved for every integer $i\in \{0,..,n-1\}$:
\be{eq:sol1}
 \begin{array}{l}
  \frac{d u_n(t)}{dt} +\Gamma^{r/ \|f\|_\infty}(C,u_n^{i+1}) \ni f_n(t) \qquad a.e.\ t\in I_i  \vsp \\
 \left\| \frac{d u_n(t)}{dt} - f_n(t) \right\| \leq \|f\|_\infty.
\end{array}
\ee

\mb
{\bf Third step: } Existence of a limit function. \\
Let $n_0$ be an integer such that Property (\ref{eq:1}) holds. First from (\ref{eq:2}) and the definition of $\Delta_n$, we deduce
that $\dot{u_n}$ is uniformly bounded by $2\|f\|_\infty$. So $(u_n)_{n\geq n_0}$ is a bounded sequence of $C([0,T],\B)$
which is uniformly Lipschitz and so it is equicontinuous. Now for each $i$ and for every $t\in I_i$, by definition we have
\be{eq:distC} d(u_n(t),C) \leq \|u_n(t)-u_n^i\| \leq \left\|u_{n}^{i+1}-u_n^i\right\| \leq h\|f\|_\infty = \frac{T}{n}\|f\|_\infty. \ee
As the set $C$ is assumed to be ball-compact and $u_n$ is bounded, we deduce that the set $\left\{ u_n(t), n\geq n_0\right\}$ is relatively compact. Then we can apply Arzela-Ascoli's Theorem to the sequence $(u_n)_n$ : there exists a subsequence, still denoted $u_n$, which converges uniformly on $[0,T]$ to a continuous function $u$. Obviously
$u(0)=u_0$. Moreover as $C$ is a closed subset, (\ref{eq:distC}) implies that the values of $u$ belong to $C$. Similarly $u$ is a Lipschitz function
and so it is absolutely continuous.

\mb
{\bf Fourth step: } The limit function $u$ is a solution of the continuous problem (\ref{eq:syst33}).\\
By the continuity of $f$, we get a pointwise convergence in $L^\infty(I,\B)$:
$$ f_n(t) \longrightarrow f(u(t)),$$
which induces the weak convergence $f_n\rightharpoonup f(u)$ in $L^\infty(I,\B)$. 
We are going to check that \be{eq:amontrer2}  \frac{du(t)}{dt}+\Gamma^{r/ \|f\|_\infty}(C,u(t)) \ni f(u(t)) \qquad a.e.\ t\in[0,T], \ee which will also imply
(\ref{eq:syst33}). \\
We have seen that the functions $\dot{u_n}$ are uniformly bounded in $L^\infty(I,\B)$. We now use our
assumption about the Banach space $\B$ with the sequence $(r'\Delta_n)_{n\geq 0}$ with $\Delta_n:=\dot{u_n}-f_n$ and $r'=r/ \|f\|_\infty$. As $\B$ is assumed to be ``$I$-smoothly weakly compact'', up to a subsequence, we may suppose without loss of generality that $(\dot{u}_n)_n$ weakly converges
in $L^\infty(I,\B)$ to a function $\omega$ such that for all $z\in L^\infty(I,\B)$ and $\phi\in L^1(I,\R)$,
\begin{align}
\lefteqn{\lim_{n\to\infty} \int_I {}_{\B^*}\langle J_p(z(t)-r'\Delta_n(t))-J_p(-r'\Delta_n(t)) , \Delta_n(t)\rangle_\B \phi(t) dt =} & & \nonumber \\
 & & \int_I  {}_{\B^*}\langle J_p(z(t)-r'\Delta(t))-J_p(-r\Delta(t)) , \Delta(t)\rangle_\B \phi(t) dt,\label{hyp2}
\end{align}
where $\Delta=\omega-f(u)$. Also $\Delta_n$ weakly converges to $\Delta$. \\
Moreover it is well-known that the weak convergence ($\dot{u_n} \rightharpoonup \omega$) implies
$$ \omega(t) = \frac{du(t)}{dt} \qquad a.e. t\in I.$$
By (\ref{eq:sol1}), we deduce that for almost every $t\in I$
$$ \left\| \dot{u}(t)-f(u(t)) \right\|_\infty \leq \|f\|_\infty.$$
Write $\tilde{u}_n(t)=u_n^{i+1}$ for $t\in I_i$ and each integer $n$. The sequence $\tilde{u}_n$ strongly converges to $u$ in $L^\infty(I,C)$. In addition for all integer $n$ and almost every $t\in I$, 
$$ \tilde{u}_n(t) \in P_C(\tilde{u}_n(t)-r'\Delta_n(t)).$$
Since Proposition \ref{prop:imp}, we deduce that this property holds for the limit functions: 
$$ u(t) \in P_C(u(t)-r'\Delta(t)), \qquad a.e. t\in I.$$
This property implies the desired one
(\ref{eq:amontrer2}) and also (\ref{eq:sys33}) which concludes the proof of the theorem. \findem

\subsection{Sweeping process in a Hilbert space for a constant set $C$.} \label{subsec:hilbert}

Here we consider a Hilbert space, denoted $\B=H$, which is a particular case of  $I$-smoothly weakly compact space (see Proposition \ref{prop:hilbert}). Before stating and proving our result, we would like to show how this assumption of a Hilbertian structure is useful. More precisely, we are going to explain how the general inequality (\ref{ineg2}) implies the ``hypomonotonicity'' property of the proximal normal cone, described by (\ref{eq:hypo}).
 Just for convenience, let us assume for this explanation that $ \|f\|_\infty=1$. For $u_0,v_0$ two initial data, we write $u$ and $v$ associated solutions (given by the previous theorem). Then (\ref{ineg2}) (used with a non constant map $\xi(t):=v(t)$ according to Remark \ref{rem:xi}) yields
\begin{align}
 \int_0^{r} \int_I
\phi(t)\left\langle J_p(u(t)-v(t)-s \Delta(t))-J_p(-s\Delta(t)),\Delta(t) \right\rangle dtds \hspace{2cm} \nonumber \\
 \hspace{6cm} \leq \frac{1}{p}\left(\int_I
\phi(t)\left\|u(t)-v(t) \right\|^p
dt\right). \label{ineg22}
\end{align}
In the case of a Hilbert space, $J_2(x)=x$ is linear (see Proposition \ref{prop:hilbert}) and so with $p=2$ we regain that
\be{eq:proxh}  \int_I
\phi(t)\left\langle u(t)-v(t),\Delta(t) \right\rangle dt \leq \frac{1}{2r}\left(\int_I
\phi(t)\left\|u(t)-v(t) \right\|^2dt\right), \ee which exactly corresponds to (\ref{eq:hypo}).
So we use $$J_2(u(t)-v(t)-s \Delta(t))-J_2(-s\Delta(t))=J_2(u(t)-v(t)).$$ We know that the linearity of $J_2$ is equivalent to a Hilbertian structure of the Banach space $\B$ (see \cite{Diestel}). 

\mb We now come to our main result. 

\begin{thm} \label{thmmh} Let $I=[0,T]$ be a bounded time-interval and $\B=H$ be a separable Hilbert space.
Let $f: H \to H$ be a bounded and Lipschitz function, $r>0$ and $C\subset H$ be a nonempty $(r,f)$-prox-regular set. Then for all $u_0\in C$, the system \be{eq:syst33h} \left\{ \begin{array}{ll}
 \dot{u}(t)+\NN(C,u(t)) \ni f(u(t)) \vsp \\
 u(0)=u_0
\end{array} \right. \ee
has one and only one absolutely continuous solution $u$, which lives in $C$. Moreover we have for almost every $t\in I$
$$\left\| \dot{u}(t)-f(u(t)) \right\|_\infty \leq \|f\|_\infty.$$
\end{thm}

\dem First we deal with the existence of solutions. \\
We will use similar arguments as for Theorem \ref{thmm}. Its proof is divided in four steps. The first, second and fourth ones did not use the ball-compactness of the set $C$ and so still hold in this case. It also remains us to develop new arguments for the third step (to prove the existence of a limit function) without requiring the ball-compactness of $C$. \\
So we refer the reader to the proof of Theorem \ref{thmm} for its steps one and two and do not recall the different notations.

\mb
{\bf New third step :} Existence of limit functions to $(u_n)_n$ and $(f_n)_n$. \\
We cannot use Arzela-Ascoli's theorem, as we do not know the relative compactness of the sets $\{u_n(t),n\geq n_0\}$. However we are going to use classical arguments (see the works cited in the introduction) to prove that $(u_n)_n$ is a Cauchy sequence in the space $L^\infty(I,H)$. We recall them in order to emphasize that these arguments, used with uniformly prox-regular sets (see \cite{Thibbv, Thibrelax} for example), still hold in the case of directional prox-regularity. \\
So for two indices $m\geq n \geq n_0$ let us consider the following function
$$\epsilon_{n,m}(t):= \left\|u_n(t)-u_m(t) \right\|^2.$$
To get an estimate of this quantity, we use Gronwall Lemma. For all $s<t \in I$, we have
$$ \int_s^t \frac{d \epsilon_{n,m}(\sigma)}{d\sigma} d\sigma = 2\int_s^t \left\langle \dot{u}_n(\sigma) - \dot{u}_m(\sigma) , u_n(\sigma)-u_m(\sigma) \right\rangle d\sigma.$$
Using the differential equation (\ref{eq:sol1}) satisfied by the discretized solutions $u_n$  and $u_m$, we have
$$ \dot{u}_n = \Delta_n + f(\tilde{u}_n) \quad \textrm{and} \quad \dot{u}_m=\Delta_m + f(\tilde{u}_m).$$
So we obtain
\begin{align}
 \int_s^t \frac{d \epsilon_{n,m}(\sigma)}{d\sigma} d\sigma = & 2 \int_s^t \left\langle \Delta_n(\sigma) - \Delta_m(\sigma) , u_n(\sigma)-u_m(\sigma) \right\rangle d\sigma \nonumber \\
  & + 2\int_s^t \left\langle f(\tilde{u}_n(\sigma)) - f(\tilde{u}_m(\sigma)) , u_n(\sigma)-u_m(\sigma) \right\rangle d\sigma. \label{eq:hy}
\end{align}
Using the Lipschitz regularity of $f$ (we denote by $L_f$ for its Lipschitz constant), we can estimate the second term of (\ref{eq:hy}) as follows~:
\begin{align*}
\int_s^t \left\langle f(\tilde{u}_n(\sigma)) - f(\tilde{u}_m(\sigma)) , u_n(\sigma)-u_m(\sigma) \right\rangle d\sigma & \\
 & \hspace{-5cm} \leq L_f \int_s^t \left\|\tilde{u}_n(\sigma) - \tilde{u}_m(\sigma) \right\| \left\| u_n(\sigma)-u_m(\sigma) \right\| d\sigma,
\end{align*}
where we have used Cauchy-Schwartz inequality.
Moreover by (\ref{eq:distC}), it can be shown that $\|\tilde{u}_n-u_n\|_\infty \leq T\|f\|_{\infty}/n $ and similarly $\|\tilde{u}_m-u_m\|_\infty \leq T\|f\|_{\infty}/m$. In using $m\geq n$, we deduce that for all $\sigma \in I$ 
\be{tilde}  \|\tilde{u}_n(\sigma) -\tilde{u}_m(\sigma)\| \leq \sqrt{\epsilon_{n,m}(\sigma)} + 2 \frac{T}{n}\|f\|_{\infty}. \ee 
Consequently, we get
\begin{align*}
\int_s^t \left\langle f(\tilde{u}_n(\sigma)) - f(\tilde{u}_m(\sigma)) , u_n(\sigma)-u_m(\sigma) \right\rangle d\sigma & \\
&  \hspace{-5cm} \leq L_f \int_s^t \left[\left\|u_n(\sigma) - u_m(\sigma) \right\| + 2\frac{T}{n}\|f\|_{\infty} \right] \left\| u_n(\sigma)-u_m(\sigma) \right\| d\sigma \\
& \hspace{-5cm} \leq L_f\int_s^t \left(\epsilon_{n,m}(\sigma) + 2\frac{T}{n}\|f\|_{\infty} \sqrt{\epsilon_{n,m}(\sigma)} \right) d\sigma.
\end{align*}
Now let us consider the first term of (\ref{eq:hy}). As $\tilde{u}_n$ and  $\tilde{u}_m$ take their values in $C$, we can apply (\ref{eq:proxh}) with $- \Delta_n\in \Gamma^{r/\|f\|_\infty}(C,\tilde{u}_n)$ and $\phi = {\bf 1}_{[s,t]}$, which gives  
$$\int_s^t \left\langle \tilde{u}_n(\sigma)-\tilde{u}_m(\sigma),\Delta_n(\sigma) \right\rangle d\sigma \leq \frac{1}{2r}\left(\int_s^t \left\|\tilde{u}_n(\sigma)-\tilde{u}_m(\sigma) \right\|^2d\sigma\right).$$
Finally by (\ref{tilde}), we deduce that
$$\int_s^t \left\langle u_n(\sigma)-u_m(\sigma),\Delta_n(\sigma) \right\rangle d\sigma \leq \frac{1}{2r} \left(\int_s^t
\left[ \sqrt{\epsilon_{n,m}(\sigma)} + 2\frac{T}{n} \|f\|_{\infty}\right]^2 d\sigma \right) + \frac{\kappa_1}{n},$$
with some constant $\kappa_1 >0$.
Similarly by symmetry, we have
$$\int_s^t \left\langle u_m(\sigma)-u_n(\sigma),\Delta_m(\sigma) \right\rangle d\sigma \leq \frac{1}{2r} \left(\int_s^t
\left[ \sqrt{\epsilon_{n,m}(\sigma)} + 2\frac{T}{n}\|f\|_{\infty}\right]^2 d\sigma\right)+ \frac{\kappa_1}{n},$$
which concludes the estimate of (\ref{eq:hy}). With the boundedness of the different sequences and of the time interval $I$, we finally have proved that
$$  \int_s^t \frac{d \epsilon_{n,m} (\sigma)}{d\sigma} d\sigma \leq 2\left(L_f + \frac{1}{r} \right)\int_s^t \epsilon_{n,m} (\sigma) d\sigma + \frac{\kappa_2}{n},$$
with some constant $\kappa_2>0$.
That holds for every $s<t$ in $I$. As $\epsilon_{n,m}(0)=0$, Gronwall's Lemma implies that
$$ \left\|u_n - u_m \right\|_\infty = \|\epsilon_{n,m} \|_\infty \leq \frac{\kappa}{n},$$
(with another constant $\kappa$),  which proves that the sequence $(u_n)_n$ is a Cauchy sequence in $L^\infty(I,H)$ and so strongly converges to a function $u$ in $L^\infty(I,H)$. This completes the ``new'' third step of the proof and we finish to show the existence of solutions in the same way as for Theorem \ref{thmm} (see the fourth part of its proof). 

\mb
{\bf Fifth step :} Uniqueness of the solutions. \\
We have seen in the explanation before the statement of the theorem, that even in the case of directional prox-regularity, the main hypomonotonicity property of $\Gamma^r(C,x)$ holds (see (\ref{eq:proxh})). So as previously, classical arguments and Gronwall's Lemma can be applied and permit to obtain the uniqueness of the solutions.
\findem

\begin{cor} \label{cor:hil} In the case of a Hilbert space, according to 
 the hypomonotonicity property of $\Gamma^r(C,x)$, all the results of Section \ref{sec1} about the equivalence to a differential equation (Proposition \ref{prop1}) and stability of solutions for (\ref{eq:sys33}), still hold with only the directional prox-regularity assumption.
\end{cor}

\subsection{Extension of our results with a moving set $C$.}

\mb In the two previous subsections, we have described two results concerning sweeping process with a constant subset $C$ under a directional prox-regularity assumption. This subsection is devoted to extend these results to the case of a moving subset $C$. Firstly, we give a generalization of Proposition \ref{prop:imp} about moving sets~:

\begin{prop} \label{prop:imp2} Let $(\B,\|\ \|)$ be a separable, reflexive and uniformly smooth Banach space. Let $C_n$ and $C:I\rightrightarrows \B$ be set-valued maps taking nonempty closed values, satisfying
\be{CnC} \sup_{t\in I} H(C_n(t),C(t)) \xrightarrow[n\to\infty]{} 0.  \ee 
We assume that for an exponent $p\in[2,\infty)$ and  a bounded sequence $(v_n)_{n\geq 0}$ of $L^\infty(I,\B)$, we can extract a subsequence $(v_{k(n)})_{n\geq 0}$ weakly converging to a point $v\in
L^\infty(I,\B)$ such that for all $z\in L^\infty(I,\B)$ and $\phi\in L^1(I,\R)$,
\begin{align}
\lefteqn{\limsup_{n\to\infty} \int_I {}_{\B^*}\langle J_p(z(t)+v_{k(n)}(t))-J_p(v_{k(n)}(t)) , v_{k(n)}(t)\rangle_\B \phi(t) dt } & & \nonumber \\
 & & \leq \int_I  {}_{\B^*}\langle J_p(z(t)+v(t))-J_p(v(t)) , v(t)\rangle_\B \phi(t) dt.\label{prop:hyp2bis}
\end{align}
Then the projection $P_{C(\cdot)}$ is weakly continuous in $L^\infty(I,\B)$ (relatively to the directions given by the sequence $(v_n)_n$) in the following sense: for all $r>0$ and for any bounded sequence $(u_n)_n$ of $L^\infty(I,\B)$ satisfying 
$$ \left\{ \begin{array}{l}
            u_n \longrightarrow u \quad \textrm{in}\ L^\infty(I,\B) \vsp \\
            u_n(t)\in P_{C_n(t)}(u_n(t)+rv_n(t))\quad \textrm{a.e.}\ t\in I,
           \end{array} \right. $$
one has for almost every $t\in I$
$$ u(t)\in P_{C(t)}(u(t)+rv(t)). $$
\end{prop}

\dem Let $\xi\in L^\infty(I,\B)$ be any map verifying $\xi(t)\in C(t)$ for all $t \in I$. Let  $\xi_n(t) \in P_{C_n(t)}(\xi(t))$ for all $t \in I$. From (\ref{CnC}), $(\xi_n)_n$ converges to $\xi$ in $L^\infty(I,\B)$ and so is bounded. The arguments of Proposition \ref{prop:imp} still hold with this non-constant map $\xi$ and permit to show that for all $\phi\in L^1(I,\R)$,
$$ \int_I \phi(t) \left[\left\| u(t)+rv(t)-\xi(t) \right\|^p-\left\|rv(t)\right\|^p\right] dt \geq 0. $$
Then we conclude as in Remark \ref{rem:xi}.
\findem

\mb The following result only requires a directional prox-regularity, but the displacement of the prox-regular set $C(\cdot)$ is supposed to be a translation.

\begin{thm} \label{thmm2} Let $\B$ be a separable, reflexive, uniformly smooth Banach space which is ``$I$-smoothly weakly compact'' for an exponent $p\in[2,\infty)$.
Let $r>0$ be a fixed real and $f: \B \to \B$ be a continuous function admitting at most a linear growth ~: there exists a constant $L>0$ such that
$$ \forall x\in \B, \qquad \|f(x)\|\leq L\left(1+\|x\|\right).$$
 Let $a\in \B$ and $C_0$ be a nonempty ``ball-compact'' and $(r,f(\cdot +ta )-a)$ prox-regular subset of $\B$ for all $t \in I$. We consider the set-valued map $C(\cdot)$ defined by 
$ \forall \, t \in I \virg C(t)=C_0+ta$.
Then for all $u_0\in C_0$, the system \be{eq:syst333} \left\{ \begin{array}{ll}
 \dot{u}(t)+\NN(C(t),u(t)) \ni f(u(t)) \vsp \\
 u(0)=u_0
\end{array} \right. \ee
has an absolutely continuous solution $u$ and for all $t\in I=[0,T]$, $u(t)\in C(t)$.
\end{thm}

{\noindent {\bf Sketch of the proof: }} The proof is essentially the same as for Theorem \ref{thmm} when $f$ is assumed to be bounded. The first step consists in defining $u_n^{i+1}=P_{C(t^{i+1})}\left[u_n^i +hf(u_n^i)\right]$ instead of (\ref{aez}), with a small enough time-step $h$ satisfying $h\|f-a\|_\infty \leq r/2$. 
As $C(t^{i+1})=C_0+t^{i+1}a$, that is equivalent to 
$$ u_n^{i+1}-t^{i+1}a = P_{C_0}\left[u_n^i-t^ia + h\left(f(u_n^i)-a\right) \right].$$
As
$$ \Delta_n(t) := \frac{1}{h} \left[u_n^{i+1}-u_n^i-hf(u_n^i) \right],$$
we have $\|\Delta_n\|_\infty \leq \|f-a \|_\infty$. Then the prox-regularity of $C_0$ in the direction $f(\cdot + t^ia)-a$ yields
$$ u_n^{i+1}-t^{i+1}a \in P_{C_0}\left[u_n^{i+1}-t^{i+1}a-\frac{r}{\|f-a\|_\infty} \Delta_n(t) \right].$$
So we get 
$$ u_n^{i+1} \in P_{C(t^{i+1})}\left[u_n^{i+1}-\frac{r}{\|f-a\|_\infty}\Delta_n(t) \right].$$
That is the key-point of the second step.
For each $i$ and $t\in I_i$, we set $C_n(t):=C(t^{i+1})$. Moreover, as for each $i$ and all $t\in I_i$, $C(t)=C_0+ta$ and $C_n(t)=C_0+(i+1)ha$, it comes~:
\be{hausd} \forall t\in I, \qquad H(C_n(t),C(t)) \leq h\|a\| =\frac{T}{n}\|a\|.\ee
Thus, the third step of the proof still holds with the following estimate~:
$$ d(u_n(t),C(t)) \leq \frac{T}{n}\left[2\|f\|_\infty + \|a\| \right].$$
The fourth step is based on Proposition \ref{prop:imp2} which can be applied because by (\ref{hausd}):
\be{h} \sup_{t\in I} H(C_n(t),C(t)) \xrightarrow[n\to\infty]{} 0. \ee 
So we obtain the existence of a solution $u$ for (\ref{eq:syst333}) satisfying for almost every $t \in I$
\be{eq:cond2} \left\| \dot{u}(t)-f(u(t)) \right\|_\infty \leq \|f-a\|_\infty. \ee
Now we explain the modifications to deal with the weaker assumption of ``linear growth'' for the perturbation $f$.
The idea (developed in \cite{Thibrelax}) is to build a sequence of maps $(u_n)_n$ which (up to a subsequence) converges uniformly to a solution of (\ref{eq:syst333}). Without loss of generality, we suppose that $4LT\leq 1$.
For every $n$, we consider a uniform subdivision $(t^i)_i$ of $I$ with a time-step $T/n$.
On $[0,t^1]$, we define $u_n$ as a solution of 
$$ \left\{ \begin{array}{ll}
 \dot{x}(t)+\NN(C(t),x(t)) \ni f(u_0) \vsp \\
 x(0)=u_0.
\end{array} \right. $$
By iterating the procedure, we define $u_n$ on $[t^i,t^{i+1}]$ as a solution of 
$$ \left\{ \begin{array}{ll}
 \dot{x}(t)+\NN(C(t),x(t)) \ni f(u_n(t^i)) \vsp \\
 x(t^i)=u_n(t^i).
\end{array} \right. $$
Then using the proof of Theorem 1 in \cite{Thibrelax}, it can be shown that
$$ \sup_{n} \max_{0\leq i<n} \|u_n(t^i)\| \leq  M,$$
for some constant $M>0$ depending on $u_0$. From (\ref{eq:cond2}), we have for almost every $t\in I_i$
$$ \left\| \dot{u}_n(t)\right\|\leq \left\| f(u_n(t^i)) - a \right\| + \| f(u_n(t^i))\| \leq 2\| f(u_n(t^i))\|+\| a\| .$$ 
The linear growth property for $f$ implies
$$  \left\| \dot{u}_n\right\|_\infty \leq 2L(1+M) + \|a\|.$$
Thus $u_n$ satisfies the following differential inclusion~:
$$  \dot{u}_n(t)+\Gamma^{\frac{r}{L(1+M)+\|a\|}}(C(t),u_n(t)) \ni f_n(t),$$
with $f_n(t)=f(u_n(t^i))$ for $t\in I_i$.
Similarly to the third and fourth steps of Theorem \ref{thmm}, we can define a limit function $u$, which will be a solution of 
$$  \dot{u}(t)+\Gamma^{\frac{r}{L(1+M)+\|a\|}}(C(t),u(t)) \ni f(u(t)),$$
according to Remark \ref{rem:xi}.
\findem

\mb We have the same extension for Theorem \ref{thmmh}, in using a similar reasoning~:

\begin{thm} \label{thmm2h} Let $\B=H$ be a separable Hilbert space.
Let $r>0$ be a fixed real and $f: H \to H$ be a Lipschitz function admitting at most a linear growth~: there exists a constant $L>0$ with
$$ \forall x\in \B, \qquad \|f(x)\|\leq L\left(1+\|x\|\right).$$
Let $a\in \B$, $r>0$ and $C_0$ be a nonempty closed $(r,f(\cdot +ta )-a)$ prox-regular subset of $H$ for all $t \in I$. We consider the set-valued map $C(\cdot)$ defined by 
$ \forall \, t \in I \virg C(t)=C_0+ta$.
Then for all $u_0\in C_0$, the system \be{eq:syst333h} \left\{ \begin{array}{ll}
 \dot{u}(t)+\NN(C(t),u(t)) \ni f(u(t)) \vsp \\
 u(0)=u_0
\end{array} \right. \ee
has one and only one absolutely continuous solution $u$ and for all $t\in I$, $u(t)\in C(t)$.
\end{thm}

\mb Without specific assumptions about the displacement of the set $C$, we have to require a uniform prox-regularity over all the directions and not only a directional one. In the framework of Banach spaces, we state the following result~:

\begin{thm} \label{thmm3} Let $\B$ be a separable, reflexive, uniformly smooth Banach space, which is ``$I$-smoothly weakly compact'' for an exponent $p\in[2,\infty)$.
Let $f: \B \to \B$ be a continuous function admitting at most a linear growth and $r>0$ be a fixed real. Let $C:t\in I \to C(t)$ be a set-valued map taking nonempty ball-compact and $r$-prox-regular values. We assume that $C(\cdot)$ moves in a Lipschitz way~: there exists a constant $k>0$ such that for all $s,t\in I$
$$ H(C(t),C(s)) \leq k|t-s|. $$
Then for all $u_0\in C(0)$, the system \be{eq:syst3333} \left\{ \begin{array}{ll}
 \dot{u}(t)+\NN(C(t),u(t)) \ni f(u(t)) \vsp \\
 u(0)=u_0
\end{array} \right. \ee
has an absolutely continuous solution $u$ and for all $t\in I$, $u(t)\in C(t)$.
\end{thm}

{\noindent {\bf Sketch of the proof: }} The proof is similar to the one of Theorem \ref{thmm2}. We only deal with the case of a bounded perturbation $f$. 
The time-step $h$ is taken in order that $h(\|f\|_\infty+k)\leq r/2$. We build the same sequence $(u_n)_n$, which satisfies~: for almost every $t\in I_i$
$$ \dot{u}_n(t) + \Gamma^{\frac{r}{\|f\|_\infty+k}}(C(t^{i+1}),u_n^{i+1}) \ni f_n(t).$$
Then we finish the proof as previously, in applying Proposition \ref{prop:imp2}~: Property (\ref{h}) is satisfied due to the Lipschitz regularity of the map $C$.
\findem


\mb We finish this article by asking the following open question: How can we get the uniqueness of sweeping process without using specific properties of a Hilbert space and how can we get around the assumption of ``ball-compactness'' of the set ? The arguments (used in Section \ref{sec1} and Subsection \ref{subsec:hilbert}) are based on Gronwall's Lemma and are specific to the Hilbert case. Mainly, the linearity of $J_2$ permits to get a very well-adapted description of the hypomonotonicity property (see (\ref{eq:proxh}). In a Banach framework, this property (called ``$J$-hypomonotonicity'') of a prox-regular set is studied by  F.~Bernard, L.~Thibault and N.~Zlateva, see \cite{BTZ,BTZ2}. However their characterizations do not allow to use Gronwall's Lemma. To obtain the uniqueness of the solutions in some Banach spaces (even in specific examples as Lebesgue spaces) seems to be a difficult problem.  We also probably need a new approach of this question.

\gb {\bf Acknowledgments :}  The two authors are indebted to
Professor Lionel Thibault for his interest about this work and valuable advices to improve this paper. Moreover they 
would like to express many thanks to him for his cordial welcome and his availability during their visit in the university of Montpellier.

\end{document}